\documentclass[12pt]{article}

\usepackage{amssymb,amsmath,amsthm, amsfonts}

\def\sqr#1#2{{\vcenter{\vbox{\hrule height.#2pt
              \hbox{\vrule width.#2pt height#1pt \kern#1pt \vrule width.#2pt}
              \hrule height.#2pt}}}}
\def\signed #1{{\unskip\nobreak\hfil\penalty50
              \hskip2em\hbox{}\nobreak\hfil#1
              \parfillskip=0pt \finalhyphendemerits=0 \par}}
\def\endpf{\signed {$\sqr69$}}
\def\3n{\negthinspace \negthinspace \negthinspace }
\def\2n{\negthinspace \negthinspace }
\def\1n{\negthinspace }

\def\={\buildrel \triangle \over =}

\def\a{\alpha}
\def\b{\beta}
\def\g{\gamma}

\def\e{\varepsilon}

\def\l{\lambda}

\def\i{\infty}

\def\O{\Omega}

\def\cD{{\cal D}}
\def\cE{{\cal E}}
\def\cF{{\cal F}}

\def\cH{{\cal H}}

\def\cL{{\cal L}}

\def\cR{{\cal R}}

\def\no{\noindent}

\def\ms{\medskip}
\def\bs{\bigskip}
\def\q{\quad}

\def\max{\mathop{\rm max}}

\def\exp{\mathop{\rm exp}}
\def\sup{\mathop{\rm sup}}

\def\cd{\cdot}

\def\inf{\hbox{\rm inf$\,$}}

\def\as{\hbox{\rm a.s.{ }}}

\def\be{\begin{equation}}
\def\bel{\begin{equation}\label}
\def\ee{\end{equation}}
\def\bt{\begin{theorem}}
\def\bcd{\begin{condition}}
\def\ecd{\end{condition}}
\def\et{\end{theorem}}
\def\bc{\begin{corollary}}
\def\ec{\end{corollary}}
\def\bde{\begin{definition}}
\def\ede{\end{definition}}
\def\bl{\begin{lemma}}
\def\el{\end{lemma}}
\def\bp{\begin{proposition}}
\def\ep{\end{proposition}}
\def\br{\begin{remark}}
\def\er{\end{remark}}
\def\ba{\begin{array}}
\def\ea{\end{array}}
\def\ed{\end{document}}

\def\square#1{\vbox{\hrule\hbox{\vrule height#1%
     \kern#1\vrule}\hrule}}
\def\rectangle#1#2{\vbox{\hrule\hbox{\vrule height#1%
     \kern#2\vrule}\hrule}}

\font\tenbb=msbm10 \font\sevenbb=msbm7 \font\fivebb=msbm5

\newfam\bbfam
\scriptscriptfont\bbfam=\fivebb \textfont\bbfam=\tenbb
\scriptfont\bbfam=\sevenbb

\newtheorem{lemma}{Lemma}[section]
\newtheorem{remark}{Remark}[section]

\newtheorem{theorem}{Theorem}[section]
\newtheorem{corollary}{Corollary}[section]

\newtheorem{definition}{Definition}[section]
\newtheorem{proposition}{Proposition}[section]
\newtheorem{condition}{Condition}[section]

\makeatletter
   
   \@addtoreset{equation}{section}
\makeatother

\begin{document}
\title{\bf Harmonic Analysis of  Stochastic Equations
and Backward Stochastic Differential Equations}

\author{Freddy Delbaen\thanks{ Department of Mathematics, Eidgen\"{o}ssische Technische Hochschule Z\"{u}rich, CH-8092 Z\"{u}rich,
Switzerland. Part of this work was done when this author was
visiting China in the years 2005, 2006, and 2007, Laboratory of
Mathematics for Nonlinear Sciences, Fudan University, whose
hospitality is greatly appreciated. Part of this work was financed
by a grant of Credit-Suisse. The paper only reflects the personal
opinion of the author. {\small\it E-mail:} {\small\tt
delbaen@math.ethz.ch}.}\q and \q Shanjian Tang\thanks{Department
of Finance and Control Sciences, School of Mathematical Sciences,
Fudan University, Shanghai 200433, China. This work is partially
supported by the NSFC under grant 10325101 (distinguished youth
foundation), the Basic Research Program of China (973 Program)
with Grant No. 2007CB814904, and the Chang Jiang Scholars Program.
Part of this work was completed when this author was visiting in
October, 2007, Department of Mathematics, Eidgen\"{o}ssische
Technische Hochschule Z\"{u}rich, whose hospitality is greatly
appreciated. {\small\it E-mail:} {\small\tt
sjtang@fudan.edu.cn}.\ms} }

\date{}

\maketitle

\begin{abstract}
The BMO martingale theory is extensively used to study  nonlinear
multi-dimensional stochastic equations (SEs) in $\cR^p$ ($p\in [1,
\infty)$) and backward stochastic differential equations (BSDEs)
in $\cR^p\times \cH^p$ ($p\in (1, \infty)$) and in
$\cR^\infty\times \overline{\cH^\infty}^{BMO}$, with the
coefficients being allowed to be unbounded. In particular, the
probabilistic version of Fefferman's inequality plays a crucial
role in the development of our theory, which seems to be new.
Several new results are consequently obtained. The particular
multi-dimensional linear case for SDEs and BSDEs are separately
investigated, and the existence and uniqueness of a solution is
connected to the property that the elementary solutions-matrix for
the associated homogeneous SDE satisfies the reverse H\"older
inequality for some suitable exponent $p\ge 1$. Finally, we
establish some relations between Kazamaki's quadratic critical
exponent $b(M)$ of a BMO martingale $M$ and the spectral radius of
the solution operator for the $M$-driven SDE, which lead to a
characterization of Kazamaki's quadratic critical exponent of BMO
martingales being infinite.


\end{abstract}

\bs

\no{\bf 2000 Mathematics Subject Classification}.  Primary 60H10,
60H20, 60H99; Secondary 60G44, 60G46.

\bs

\no{\bf Key Words}. BMO martingales, stochastic equations,
backward stochastic differential equations, Fefferman's
inequality, reverse H\"{o}lder inequalities, unbounded
coefficients.


\section{Preliminaries}

Let $T>0$. Let $(\O,\cF,\{\cF_t\}_{t\ge 0},P)$  be a complete
filtered probability space on which a one-dimensional standard
Brownian motion $\{W(t)\}_{t\ge0}$ is defined such that
$\{\cF_t\}_{t\ge 0}$ is the natural filtration generated by
$W(\cd)$, augmented by all the $P$-null sets in $\cF$. Let $H$ be
a Banach space. We denote by $\cL_{\cF}^p(0,T;H)$ ($p\ge 1$) the
Banach space consisting of all $H$-valued $\{\cF_t\}_{t\ge
0}$-optional processes $X(\cd)$ such that
$E(|X(\cd)|_{L^p(0,T;H)}^2)<\i$, with the canonical norm; by
$\cL_{\cF}^\infty (0,T;H)$ the Banach space consisting of all
$H$-valued $\{\cF_t\}_{t\ge 0}$-optional bounded processes; and by
$\cL_{\cF}^2(\O;C([0,T];H))$ the Banach space consisting of all
$H$-valued $\{\cF_t\}_{t\ge 0}$-adapted continuous processes $X$
such that $E(|X|_{C([0,T];H)}^2)<\i$, with the canonical norm.

\bde Let $p\in [1,\infty)$. The space $\cR^p$ is the space of all
continuous adapted processes $Y$ such that \be \|
Y||_{\cR^p}:=\|Y^*_T\|_{L^p} \quad \hbox{\rm with }
Y^*_T:=\max_{0\le t\le T}|Y_t|\ee is finite. $\cH^p$ is the Banach
space of continuous $\{\cF_t, 0\le t\le T\}$-adapted local
martingales such that \be \|Y\|_{\cH^p}:=\|\langle
Y\rangle^{1/2}_T\|_{L^p}\ee is finite. $\langle Y\rangle$ denotes
the quadratic variation process of a semi-martingale, and $\langle
X,Y\rangle$ denotes the covariance process between the two
semi-martingales $X$ and $Y$. \ede

Let $M$ be a continuous martingale. Define \be a(M):=\sup \{a\ge
0: \sup_{\tau} \left\|E\left[\exp
\left(a|M_\infty-M_\tau|\right)|\cF_\tau\right]\right\|_{L^\infty}<\infty\}\ee
and \be b(M):=\sup \left\{b\ge 0: \sup_{\tau}
\left\|E\left[\exp{\left({1\over2}b^2\left(\langle
M\rangle_\infty-\langle M\rangle_\tau \right)\right)\biggm
|\cF_\tau}\right]\right\|_{L^\infty}<\infty\right\}. \ee In both
expressions, $\tau$ is an arbitrary stopping time.

\bde Let $Y=(Y_t)_{0\le t\le T}$ be a uniformly integrable
martingale. Then Y is said to {\bf belong to BMO} if there is a
constant $C>0$ such that for every stopping time $\tau$\be
E\left[|Y_T-Y_\tau|^p|\cF_\tau\right] \le C \quad P\hbox{\rm
-}\as.\ee  This definition is independent of $p$. Usually we
define $\|Y\|_{BMO}$ as the smallest constant $c$ such that for
all stopping time $\tau$, \be
E\left[|Y_T-Y_\tau|^2|\cF_\tau\right] \le c^2 \quad P\hbox{\rm
-}\as. \ee \ede


 \bde The nonzero-valued process $L$
is said to  satisfy the reverse H\"{o}lder inequality under $P$,
denoted by $R_p(P)$, where $p\in [1, +\infty]$, if there is a
constant $C>0$ such that for every  stopping time $\tau$, we have
\be E\left[\left|{L_T\over L_\tau}\right|^p\biggm|
\cF_\tau\right]\le C. \ee For $p=+\infty$, we require that
${L_T\over L_\tau}$ is essentially bounded by $C$
 (see Kazamaki~\cite[Definition 3.1. ]{Kaz}).\ede

\bl \label{KW} (Kunita-Watanabe inequality) Let $X$ and $Y$ be two
semi-martingales, and let $H$ and $K$ be two measurable processes.
Then, we have almost surely \be \int_0^\infty |H_s||K_s||d
[X,Y]_s|\le \left(\int_0^\infty H_s^2d
[X,X]_s\right)^{1/2}\left(\int_0^\infty K_s^2d
[Y,Y]_s\right)^{1/2}.\ee More generally, for $p\in [1,\infty)$, we
have \be\label{KWH} \int_0^\infty |H_s||K_s||d [X,Y]_s|\le
\left(\int_0^\infty H_s^pd [X,X]_s\right)^{1/p}\left(\int_0^\infty
K_s^qd [Y,Y]_s\right)^{1/q}\ee with ${1\over p}+{1\over q}=1.$\el

\bl (Fefferman's inequality) If $X\in \cH_1$ and $Y\in BMO$, then
\be E\left[\int_0^T|d\langle X,Y\rangle_s|\right]\le \sqrt{2}
\|X\|_{\cH^1}\|Y\|_{BMO}.\ee\el

About the expression of the duality between $\cH^1$ and BMO space,
we have (see Kazamaki~\cite[Theorem 2.7, page 38]{Kaz}):

\bl \label{1-infty duality} Let $X$ be a continuous local
martingale. Then, we have
\be \ba{rcl} \|X\|_{\cH^1}&\le& \sup\{E\left[\langle X,Y\rangle_\infty\right]: \|Y\|_{BMO}\le 1\},\\
 \|X\|_{BMO}&\le& \sup\{E\left[\langle X,Y\rangle_\infty\right]: \|Y\|_{\cH^1}\le
 1\}.
\ea\ee\el

From Fefferman's inequality, we can show the following lemma.

\bl \label{Emery} Let $p\in [1, \infty)$. Assume that $X\in \cR^p$
and $M\in BMO$. Then, $X\circ M\in \cH^p$. Moreover, we have the
following estimate \be\|X\circ M\|_{\cH^p} \le
\sqrt{2}\|X\|_{\cR^p}\|M\|_{BMO}. \ee for $p\in (1, \infty)$ and
\be\|X\circ M\|_{\cH^1} \le \|X\|_{\cR^1}\|M\|_{BMO} \ee
(corresponding to the case of $p=1$).\el

{\bf Proof of Lemma~\ref{Emery}. } (i) The case $p\in (1,
\infty)$. Take any $N\in \cH^q$. We have \be\ba{rcl}&&
E\left[|\langle X\circ M, N\rangle_\infty|\right]\le
E\left[|\langle X\circ
N, M\rangle_\infty|\right]\\
&\le & \sqrt{2}\|X\circ N\|_{\cH^1}\|M\|_{BMO} \quad \hbox{\rm (using Fefferman's inequality)}\\
&\le&\sqrt{2}\|X\|_{\cR^p}\|N\|_{\cH^q}\|M\|_{BMO}. \quad
\hbox{\rm (using H\"older's inequality)}\ea \ee

(ii) The case $p=1$. We have \be \ba{rcl} \displaystyle
\int_0^\infty X_s^2d \langle M \rangle_s&\le &\displaystyle
X^*_\infty\int_0^\infty|X_s|d\langle M\rangle_s\\[0.3cm]
&\le & \displaystyle X^*_\infty\int_0^\infty X_s^*d\langle M\rangle_s\\[0.3cm]
&\le &\displaystyle  X^*_\infty\left (X^*_\infty\langle
M\rangle_\infty-\int_0^\infty \langle M\rangle_s
dX_s^*\right)\\[0.3cm]
&\le &\displaystyle  X^*_\infty\left( \int_0^\infty \left(\langle
M\rangle_\infty-\langle M\rangle_s\right) dX_s^*\right). \ea\ee
Therefore, \be\ba{rcl} \displaystyle E\left[\left(\int_0^\infty
X_s^2d\langle M \rangle_s\right)^{1/2}\right]&\le &\displaystyle
E\left[\left(X_\infty^*\int_0^\infty\left(\langle
M\rangle_\infty-\langle
M\rangle_s \right)dX_s^*\right)^{1/2}\right]\\
&\le&\displaystyle
\left\{E\left[X_\infty^*\right]\right\}^{1/2}\left\{E\left[\int_0^\infty(\langle
M\rangle_\infty-\langle M \rangle_t)\, dX_t^*\right]\right\}^{1/2}\\
&\le&\displaystyle  \|X\|_{\cR^1}^{1/2}\left\{E\left[\int_0^\infty
E\left[\left(\langle
M\rangle_\infty-\langle M \rangle_t\right)|\cF_t\right]\, dX_t^*\right]\right\}^{1/2}\\
&\le & \|X\|_{\cR^1}^{1/2}
\|M\|_{BMO}\left\{E\left[X_\infty^*\right]\right\}^{1/2}\le
\|X\|_{\cR^1}\|M\|_{BMO}. \ea \ee The proof is complete.
\endpf

For the case of $X\in \cH^p$ ($\subset \cR^p$), the first
assertion in Lemma~\ref{Emery} is included in Ba\~nuelos and
Bennett~\cite[Theorem 1.1 (i), page 1227]{BanBenn}. The following
lemma  is obvious from the definition of $BMO$ norm, see
Ba\~nuelos and Bennett~\cite[Theorem 1.1 (ii), page
1227]{BanBenn}.

\bl  If $X\in \cR^\infty$ and $M\in BMO$, then $X\circ M\in BMO$
and $\|X\circ M\|_{BMO}\le \|X\|_{\cR^\infty} \|M\|_{BMO}.$ \el

\bl \label{Lp} Let $p\in [1, \infty)$. Assume that $X\in \cH^p$
and $M\in BMO$. Then, $\langle X, M\rangle_\infty \in L^p$.
Moreover, we have the following estimate \be\|\langle X,
M\rangle_\infty\|_{L^p} \le \sqrt{2}p \|X\|_{\cH^p}\|M\|_{BMO}.\ee
\el

The first assertion in Lemma~\ref{Lp} can be found in Ba\~nuelos
and Bennett~\cite[Theorem 1.1 (iii), page 1227]{BanBenn}. For
convenience of the reader, we give a full proof.

 {\bf
Proof of Lemma~\ref{Lp}. } For the case $p=1$, noting that \be
|\langle X,M\rangle_\infty|\le \int_0^\infty |d\langle
X,M\rangle|,\ee
 it is immediate from Fefferman's inequality to get
the desired results. In what follows, we consider the case  $p\in
(1, \infty)$. Then, $q\in (1, \infty)$. Take any $\xi\in L^q$.
Write $Y_t:=E[\xi|\cF_t]$ for $t\in [0,\infty]$. We have
$Y_\infty=\xi$ and \be\ba{rcl}&&\displaystyle  E\left[\langle X,
M\rangle_\infty \xi\right]= E\left[\int_0^\infty Y_sd\langle X,
M\rangle_s \right]= E\left[\int_0^\infty d\langle X, Y\circ
M\rangle_s \right]\\
&\le & \|X\|_{\cH^p}\|Y\circ M\|_{\cH^q} \\[0.2cm]
&& \quad \quad \hbox{\rm (using both Kunita-Watanabe inequality
and H\"older's
inequality)}\\[0.2cm]
&\le & \sqrt{2}\|X\|_{\cH^p}\|M\|_{BMO} \|Y\|_{\cR^q}\quad \hbox{\rm (using Lemma~\ref{Emery})}\\
&\le&\sqrt{2}p\|X\|_{\cH^p}\|\xi\|_{L^q}\|M\|_{BMO}. \quad
\hbox{\rm (using Doob's inequality)}\ea \ee
\endpf

\bde An integrable random variable $\xi$ is said to be in $BMO$ if
the local martingale $\{E[\xi |\cF_t], t\in [0,T]\}\in BMO.$  \ede

\bl \label{Linfty}Let $X\in BMO$ and $M\in BMO$. Then, $\langle
X,M\rangle_\infty\in BMO$. Moreover, $\|\langle
X,M\rangle_\infty\|_{BMO}\le \sqrt{2}\|X\|_{BMO}\|M\|_{BMO}.$\el

{\bf Proof. } Take $Y\in \cH^1$. We have \be
\ba{rcl}&&\displaystyle  \left|E\left[Y\langle X,M\rangle
\right]\right| =\left|E\left[\int_0^\infty Y_s\, d\langle X,M
\rangle_s\right ]\right|\\
&=&\displaystyle  \left|E\left[\langle Y\circ
X,M\rangle\right]\right |\le \sqrt{2} \|Y\circ
X\|_{\cH^1}\|M\|_{BMO}\quad \hbox{\rm (Fefferman's inequality) } \\
&\le &\displaystyle  \sqrt{2} \|Y\|_{\cH^1}\|X\|_{BMO}\|M\|_{BMO}.
\quad \hbox{\rm (Lemma~\ref{Emery}) } \ea\ee Using
Lemma~\ref{1-infty duality}, we have the desired results.
\endpf

The following fundamental Burkh\"older-Davis-Gundy (abbreviated as
BDG) inequality will be frequently used in our paper: for any
$p\in (0,\infty)$, there are two universal positive constants
$c_p$ and $C_p$ such that for any local continuous martingale $M$
with $M_0=0$, we have \be C_p^{-p} E\left[\langle M\rangle_T
^{p/2}\right] \le E\left[\left(M^*_T\right)^p\right] \le c_p^{-p}
E\left[\langle M\rangle_T ^{p/2}\right],\ee or in a different
form, \be C_p^{-1} \|M\|_{\cH^p} \le \|M\|_{\cR^p} \le
c_p^{-1}\|M\|_{\cH^p} .\ee See Yor~\cite[page 100]{Yor}.

The following definition is based on that of
Emery~\cite{Emery1,Emery2} (see also Protter~\cite[page
248]{Protter}).

\bde Let $M\in BMO$ and $\varepsilon>0$. A finite sequence of
stopping times $0=T_0\le T_1\le \cdots\le T_k$ is said to
$\varepsilon$-slice $M$ if $M=M^{T_k}$ and
$|(M-M^{T_i})^{T_{i+1}}|_{BMO}\le \varepsilon$, for
$i=0,1,\cdots,k-1$. If such a sequence of stopping times exists,
we say that $M$ is $\varepsilon$-sliceable in $BMO$. \ede

\bde $M$ is called sliceable in $BMO$ if for $\forall \e>0$, $M$
is $\varepsilon$-sliceable in $BMO$, i.e., there are a positive
integer $N$ and a finite increasing sequence of stopping times
$\{T_i, i=1,2, \ldots, N. \}$ with $T_0=0$ and $T_{N+1}=\infty$
such that ${}^{T_n}M^{T_{n+1}}:=M^{T_{n+1}}-M^{T_n}$ satisfies \be
\|{}^{T_n}M^{T_{n+1}}\|_{BMO}\le \e. \ee This is equivalent to $M
\in \overline{\cH^\infty}^{BMO}$ by Schachermayer's
result~\cite{Schacher}. \ede

For more knowledge on local martingales and semi-martingales, the
reader is referred to, among others, the following books:
Dellacherie and Meyer~\cite{DellMeyer}, He, Wang, and
Yan~\cite{HeWangYan}, Kazamaki~\cite{Kaz}, and
Protter~\cite{Protter}.

Throughout the rest of the paper, $N_1,N_2,$ and $M$ are supposed
to be continuous local martingales on the time interval $[0,T]$,
being equal to zero at time $t=0$.

\vspace{1cm} Since It\^o's initial works~\cite{Ito1,Ito2,Ito3},
stochastic differential equations (abbreviated hereafter as SDEs)
driven by general semimartingales, instead of just Brownian
motion, have been studied by Dol\'eans-Dade~\cite{Dol},
Dol\'eans-Dade and Meyer~\cite{DolMeyer},
Protter~\cite{Protter1,Protter}, and Emery~\cite{Emery2,Emery2}
among others. The theory of existence and uniqueness on SEs driven
by general semi-martingales is already quite general.  However,
the rather general result presented in the literature  is
concerned with existence and uniqueness in a very large space like
$\cup_{p\ge 1}\cH^p$. In this subsection, we present some new
sufficient conditions on existence and uniqueness of solutions in
$\cH^p$ for some fixed $p\in [1,\infty)$. These conditions  are
more general than those presented in Protter~\cite{Protter},
allowing the coefficients to be unbounded. We make best use of the
deep property of Fefferman's inequality on BMO martingales, which
seems to be new in the study of SEs.

Similar situations also exist for the research into BSDEs. Since
Bismut's initial works~\cite{Bis1,Bis2,Bis3} and Pardoux and
Peng's seminal paper~\cite{PP}, BSDEs driven by general local
martingales in the space $\cR^p\times \cH^p$ for general $p\in
(1,\infty)$ instead of just $p=2$, have been studied by
Buckdahn~\cite{Buck} (with the restriction that $p\in [2,
\infty)$) and El Karoui, Peng and Quenez~\cite{ElPengQuenez} (the
underlying driving martingale is assumed to be a Brownian motion)
among others. In  El Karoui, Peng and Quenez~\cite{ElPengQuenez},
the coefficients of BSDEs are  restricted to be uniformly
Lipschitz in the unknown variables.  The existence results in the
space $\cR^p\times \cH^p$ for some $p\in [2,\infty)$ existing in
Buckdahn~\cite{Buck} requires---though the coefficients of BSDEs
are allowed to be unbounded
---that the data $(\xi, J)$ (see BSDE~(\ref{bsde-basic}) below) lie
in a space $\cR^{p+\epsilon}$ for some $\epsilon>0$,  a stronger
integrability. Roughly speaking, the integrability of the adapted
solution of BSDEs is less than that of the data in
Buckdahn~\cite{Buck}. Note that BSDEs with unbounded coefficients
have also been studied by El Karoui and Huang~\cite{ElHuang}, but
requiring that both the solution and the data lie in the square
integrable space which is weighted in relevance to the
coefficients. In this paper, the BMO martingale theory, in
particular Fefferman's inequality on BMO martingales, is applied
to study BSDEs with unbounded coefficients. New existence results
are proved where the adapted solutions of BSDEs---even though the
coefficients are unbounded---have the same integrability index $p$
to the underlying data $(\xi, J)$ for $p\in (1,\infty)$. The
critical case of $p=+\infty$ is also discussed, and some
interesting results are obtained.

It seems to be necessary to mention some applications of BMO
martingales in the study of BSDEs. Bismut~\cite{Bis3} has already
used some properties of BMO matingales when he discussed the
existence and uniqueness of adapted solutions of backward
stochastic Riccati equation in some particular case. He chose the
BMO space for the second unknown variable. In the work of Delbaen
et al. ~\cite{DMSSS1,DMSSS2} on hedging contingent claims in
mathematical finance, BMO martingales are connected to some
closedness in some suitable Banach spaces of the set of attainable
claims for the agent's wealth equation, which is essentially  a
problem of existence and uniqueness of a linear BSDE, but with
unbounded coefficients. In the conference on mathematical finance,
held in Konstanz in the year of 2000, the role of BMO martingales
received a special emphasis in the study of backward stochastic
Riccati equation and related linear quadratic stochastic optimal
control problems. See Kohlmann and Tang~\cite{KT1,KT2,KT3}. In
particular in Kohlmann and Tang~\cite{KT3}, the second component
of the adapted solution pair for a general  backward stochastic
Riccati equation---which is a multi-dimensional BSDE with the
generator being a quadratic form of the second unknown variable
---is shown to be a BMO martingale. Later, such kind of  results
are widely obtained and used, among others, by Hu, Imkeller, and
M\"uller~\cite{HIM}, Hu and Zhou~\cite{HuZhou1}, Barrieu and El
Karoui~\cite{BarrEl}, Briand and Hu~\cite{BrHu1,BrHu2}, and Hu et
al.~\cite{HMPY}.

The rest of the paper consists of three sections, and is organized
as follows.

Section 2 consists of three subsections. In Subsection 2.1, a
rather general nonlinear multi-dimensional SE~(\ref{nsdeG}) driven
by semimartingales with unbounded coefficients is discussed, and a
new existence result in $(\cR^p)^n$ ($p\in [1, \infty)$) is proved
under some suitable sliceability in the BMO space of the
coefficients, which is stated in Theorem~\ref{nse}. In Subsection
2.2, a rather general nonlinear multi-dimensional
BSDE~(\ref{nbsde}) driven by a continuous local martingale with
unbounded coefficients is discussed, and a new  existence result
in $(\cR^p)^n\times (\cH^p)^{(2n)}$ ($p\in (1, \infty)$) is proved
under some suitable sliceability in the BMO space of the
coefficients, which is stated in Theorem~\ref{bsde-basic}. For the
critical case of $p=\infty$, a new existence result in
$\cap_{p>1}(\cR^p)^n\times (BMO)^{(2n)}$ is also obtained, but for
a less general BSDE~(\ref{nbsdeinfty}), and it is stated in
Theorem~\ref{bsde-basic-infty}. In Subsection 2.3, we give a
sufficient condition on the suitable sliceability in the BMO space
of the coefficients required in Theorems~\ref{nse},
\ref{bsde-basic}, and \ref{bsde-basic-infty}. They are stated in
Theorems~\ref{nsesliceable}, \ref{bsde-basic-sliceable}, and
\ref{bsde-basic-infty-sliceable}, respectively. Moreover, when the
data~$(\xi,J)\in (L^\infty(\cF_T))^n\times (\cR^\infty)^n$, a new
existence result in $(\cR^\infty)^n\times
\left(\overline{\cH^\infty}^{BMO}\right)^{(2n)}$ is proved for the
rather general nonlinear multi-dimensional BSDE~(\ref{nbsde}) with
a nice application of Fefferman's inequality, the John-Nirenberg
inequality, and the Garnett-Jones's Theorem, and it is stated in
Theorem~\ref{nbsdebounded}.

Section 3 is concerned with the linear BSDEs and SDEs with
unbounded coefficients.  The existence and uniqueness of the
solution is connected to some reverse H\"older inequality
property. It consists of two subsections. Subsection 3.1 is
concerned with linear BSDEs with unbounded coefficients, while
Subsection 3.2 is concerned with linear SDEs with unbounded
coefficients.

Finally, in Section 4, the solution operator $\phi$ from $\cH^p$
to $\cH^p$ of the one-dimensional SDE driven by a BMO martingale
$M$ receives a special consideration, whose spectral radius is
estimated in terms of the Kazamaki's quadratic critical exponent
$b(M)$ for the underlying BMO martingale $M$. This estimation
leads to a characterization of $b(M)=\infty$.

\section{The nonlinear multi-dimensional case}

\subsection {Unbounded SEs}

Let $\cD$ denote the space of $\{\cF_t, 0\le t\}$-adapted
c\`adl\`ag processes, and $\cD^n$ the space of $n$-dimensional
vector processes whose components are in $\cD$.

Consider the following nonlinear SEs: \be
X_t=J(t)+\int_0^tf(s,X)\, d\langle N_1,N_2\rangle_s+\int_0^tg(s,
X)\, dM_s, \quad t\in[0,T].\label{nsdeG}\ee Here,  $J\in
(\cR^p)^n, f$ and $g$ denote $R^n$-valued functionals defined on
$\Omega\times [0,T]\times \cD^n.$

\bt \label{nse} Let $p\in [1, \infty)$. Assume that

(i) There are two $\{\cF_t, 0\le t\le T\}$-adapted processes
$\a(\cdot)$ and $\b(\cdot)$ such that \be f(t,0)=0; \quad
|f(t,x_1)-f(t,x_2)|\le \a (t)(x_1-x_2)^*(t), \quad x_1,x_2\in
\cD^n\ee and \be g(t,0)=0; \quad |g(t,x_1)-g(t,x_2)|\le \b
(t)(x_1-x_2)^*(t), \quad x_1,x_2\in \cD^n.\ee

(ii) The martingale  $\a \circ N_1\in BMO$. The martingale $
N_2\in BMO$ is $\varepsilon_1$-sliceable in the space $BMO$ and
the martingale $\b\circ M\in BMO$ is $\varepsilon_2$-sliceable in
the space $BMO$. Let \be \label{contractionCON} \rho_1:=
2p\varepsilon_1 |\a\circ N_1|_{BMO}+\sqrt{2}\varepsilon_2
C_p<1.\ee

Then for any $J\in (\cR^p)^n$, there is unique solution $X\in
(\cR^p)^n$ to equation~(\ref{nsdeG}). Furthermore, there is a
constant $K_p$, which is independent of $J$, such that \be
\label{aprioriinqual} \|X\|_{\cR^p}\le K_p \|J\|_{\cR^p}.\ee If
$J\in (\cR^p)^n$ is a semi-martingale, then so is the solution.
\et

{\bf Proof. } We shall use the contraction mapping principle to
look for a fix-point. For this purpose, consider the following map
$I$ in the Banach space  $(\cR^p)^n$ : \be\label{mapdef}
I(X)_t:=J(t)+\int_0^tf(s,X)\, d\langle
N_1,N_2\rangle_s+\int_0^tg(s,X)\, dM_s,\quad t\in [0,T].\ee We
have \be \label{Est1} \ba{rcl} &&\displaystyle E\left[\max_{0\le
t\le T}\left|\int_0^tf(s,X)\, d\langle
N_1,N_2\rangle_s\right|^p\ \right]\\[0.3cm]
&\le&\displaystyle E\left[\left|\int_0^T|f(s,X)|\, |d\langle N_1,N_2\rangle|_s\right|^p \ \right]\\[0.3cm]
&\le&\displaystyle E\left[\left|\int_0^T\a_s X_s^*\, |d\langle N_1,N_2\rangle|_s\right|^p\ \right] \\[0.3cm]
&\le&\displaystyle  (\sqrt{2}p)^p\|\a\circ N_1\|_{BMO}^p\|X^*\circ
N_2\|_{(\cH^p)^n}^p \quad  \hbox{\rm (using Lemma~\ref{Lp})}\\
&\le & \displaystyle   (2p)^p\left\|\a\circ
N_1\right\|_{BMO}^p\|X\|_{(\cR^p)^n}^p\|N_2\|^p_{BMO} \quad
\hbox{\rm (using Lemma~\ref{Emery})}\ea\ee and \be
\label{Est2}\ba{rcl}
&&\displaystyle E\left[\max_{0\le t\le T}\left|\int_0^tg(s,X)\, dM_s\right|^p\ \right]\\[0.3cm]
&\le&\displaystyle  C_p^pE\left[\left|\int_0^T|g(s,X)|^2\,
d\langle M\rangle_s\right|^{p/2}\right]\quad \hbox{\rm (from the
BDG
inequality)} \\[0.3cm]
&\le & \displaystyle  C_p^pE\left[\left|\int_0^T\b_s^2(X_s^*)^2\,
d\langle
M\rangle_s\right|^{p/2}\right] \quad \hbox{\rm (from the Lipschitz assumption on $g$ )} \\[0.3cm]
&= &\displaystyle C_p^p \|(\b X^*)\circ M\|^p_{(\cH^p)^n}\\
&= &\displaystyle C_p^p \| X^*\circ (\b\circ M) \|^p_{(\cH^p)^n}\\
&\le &\displaystyle  (\sqrt{2} C_p)^p\|\b\circ
M\|^p_{BMO}\|X\|^p_{(\cR^p)^n}. \quad \hbox{\rm (using
Lemma~\ref{Emery})}\ea\ee
 Therefore, $I(X)\in (\cR^p)^n$ for $X\in
(\cR^p)^n$.

For $X^1,X^2\in (\cR^p)^n$, proceeding similarly to the above
arguments, we have \be \ba{rcl} &&\displaystyle E\left[\max_{0\le
t\le
T}\left|\int_0^t\left[f(s,X^1)-f(s,X^2)\right]\, d\langle N_1,N_2\rangle_s\right|^p\ \right]\\[0.3cm]
&\le & \displaystyle  (2p)^p\|\a\circ
N_1\|_{BMO}^p\|X^1-X^2\|_{(\cR^p)^n}^p\|N_2\|^p_{BMO} \ea \ee and
\be \ba{rcl}
&&\displaystyle E\left[\max_{0\le t\le T}\left|\int_0^t\left[g(s,X^1)-g(s,X^2)\right]\, dM_s\right|^p\ \right]\\[0.3cm]
&\le &\displaystyle (\sqrt{2} C_p)^p\|\b\circ
M\|^p_{BMO}\left\|X^1-X^2\right\|^p_{(\cR^p)^n}.\ea\ee Therefore,
we have
\be\label{estimate} \ba{rcl}&&\displaystyle \left\|I(X^1)-I(X^2)\right\|_{(\cR^p)^n} \\[0.3cm]
&\le&\displaystyle  \left[\left(\sqrt{2}C_p\right)^p\|\b\circ
M\|_{BMO}^p+(2 p)^p\|\a\circ
N_1\|_{BMO}^p\|N_2\|^p_{BMO}\right]^{1/p}\left\|X^1-X^2\right\|_{(\cR^p)^n}
\\[0.3cm]
&\le&\displaystyle   \left[\sqrt{2}C_p\|\b\circ M\|_{BMO}+
2p\|\a\circ
N_1\|_{BMO}\|N_2\|_{BMO}\right]\left\|X^1-X^2\right\|_{(\cR^p)^n}.
\ea \ee

Since the martingale $N_2\in BMO$ is $\varepsilon_1$-sliceable and
$\b\circ M\in BMO$ is $\varepsilon_2$-sliceable,  there is a
finite sequence of stopping times $\{T_i, i=1,2,\cdots,\widetilde
I \}$ such that the following are satisfied:

(i) $0=T_0\le T_1 \le T_2\le \cdots\le T_{\widetilde I}\le
T_{\widetilde I+1}=T$;

(ii) $|N_{2i}|_{BMO}\le \varepsilon_1$, and $|\b\circ
M_i|_{BMO}\le \varepsilon_2$ where
$N_{2i}:=N_2^{T_{i+1}}-N_2^{T_i}$ and $M_i:=M^{T_{i+1}}-M^{T_i}$
are defined on $[T_i,T_{i+1}]$. Since $|\a\circ N_{1i}|_{BMO}\le
|\a\circ N_1|_{BMO}$, we have \be
\label{contraction}\rho_{1i}:=2p\varepsilon_1 |\a\circ
N_{1i}|_{BMO}+\sqrt{2}\varepsilon_2 C_p\le\rho_1 \ee with
$N_{1i}:=N_1^{T_{i+1}}-N_1^{T_i}$ for $i=0,1, 2,\cdots,{\widetilde
I}$. Set $\cR^p_i:=\cR^p[T_i,T_{i+1}]$. Set $X_{-1}:=0$. Consider
the map $I_i:(\cR^p_i)^n\to (\cR^p_i)^n$, defined by \be
I_i(X)_t:=J_i(t)+\int_{T_i}^tf(s,X_s)\, d\langle
N_{1i},N_{2i}\rangle_s+\int_{T_i}^tg(s,X_s)\, dM_{is},\quad t\in
[T_i,T_{i+1}],\ee where
$J_i(\cdot):=J^{T_{i+1}}-J(T_i)+X_{i-1}(T_i)$ is defined on
$[T_i,T_{i+1}]$.

Similar to the derivation of
inequality~(\ref{estimate}), we have \be\ba{rcl}&&\displaystyle \left\|I_i(X^1)-I_i(X^2)\right\|_{(\cR^p_i)^n} \\
&\le&\displaystyle   \left[\sqrt{2}C_p\|\b\circ M_i\|_{BMO}+
2p\|\a\circ
N_{1i}\|_{BMO}\|N_{2i}\|_{BMO}\right]\left\|X^1-X^2\right\|_{(\cR_i^p)^n}\\[0.3cm]
&\le & \displaystyle   \left[\sqrt{2}C_p\varepsilon_2+
2p\varepsilon_1\|\a\circ
N_{1i}\|_{BMO}\right]\left\|X^1-X^2\right\|_{(\cR_i^p)^n}\\[0.3cm]
&= & \displaystyle \rho_{1i}
\left\|X^1-X^2\right\|_{(\cR_i^p)^n}\le \rho_1
\left\|X^1-X^2\right\|_{(\cR_i^p)^n}\ea \ee for any $X^1,X^2\in
(\cR^p_i)^n$. In view of the second assumption of the theorem,
 we see that the map $I_i$ is a contraction map, and satisfies the
following estimate: \be \|I(X)\|_{(\cR_i^p)^n} \le \rho_1
\|X\|_{(\cR_i^p)^n}+\|J\|_{(\cR_i^p)^n} \ee for any $X\in
(\cR^p_i)^n$. Therefore, in an inductive way, we show that the
following stochastic equation \be X_t=J_i(t)+\int_{T_i}^tf(s,X)\,
d\langle N_{1i},N_{2i}\rangle_s+\int_{T_i}^tg(s,X)\, dM_{is},\quad
t\in [T_i,T_{i+1}]\ee has a unique solution $X_i(\cdot)$ in
$(\cR_i^p)^n$ for $i=0,1,\cdots, \widetilde I$. Moreover, we have
\be\label{Est3} \|X_i\|_{(\cR^p_i)^n}\le
(1-\rho_1)^{-1}\|J_i\|_{(\cR^p_i)^n}. \ee

Then, the process \be X(t):=\sum_{i=0}^{\widetilde
I}X_i(t)\chi_{[T_i, T_{i+1})}(t), \quad t\in [0,T]\ee lies in
$(\cR^p)^n$ and is the unique solution to equation~(\ref{nsdeG}).
The desired a priori estimate (\ref{aprioriinqual}) is immediate
from the assumption (\ref{contractionCON}) and the  inequality
(\ref{Est3}). The last assertion of the theorem is obvious.
\endpf

\subsection{Unbounded BSDEs}

Consider the following nonlinear BSDEs: \be\left\{\ba{rcl}
Y_t&=&\displaystyle \xi+J_T-J_t+\int_t^Tf(s, Y_s)\, d\langle
N_1,N_2\rangle_s+\int_t^Tg(s,Y_s,Z_s)\, d\langle
M\rangle_s\\[0.3cm]
&&\displaystyle -\int_t^TZ_s\, dM_s-\int_t^T dM^\bot_s, \quad t\in
[0,T];\quad \langle M, M^\bot\rangle=0. \ea
\right.\label{nbsde}\ee Here, $\xi$ is an $R^n$-valued
$\cF_T$-measurable  random variable, $J$ is an $R^n$-valued
optional continuous process, and the $R^n$-valued random fields
$f$ and $g$ are defined on $\Omega\times [0,T]\times R^n$ and
$\Omega\times [0,T]\times R^n\times R^n$, respectively.

\bt\label{bsde-basic} Let $p\in (1, \infty)$ and $q$ be the
conjugate number. Assume that

(i) There are three $\{\cF_t, 0\le t\le T\}$-adapted processes
$\a(\cdot), \b(\cdot)$ and $\g(\cdot)$ such that \be f(\cdot,0)=0;
\quad |f(t,y_1)-f(t,y_2)|\le \a (t)|y_1-y_2|\ee for $ y_1,y_2\in
R^n$ and \be g(\cdot,0,0)=0; \quad |g(t,y_1,z_1)-g(t,y_1,z_2)|\le
\b (t)|y_1-y_2|+ \g (t)|z_1-z_2|\ee for $y_1,y_2, z_1,z_2\in R^n.$

(ii) The martingale  $\a \circ N_1\in BMO$. The martingale $
N_2\in BMO$ is $\varepsilon_1$-sliceable in the space $BMO$, the
martingale $ \sqrt{\b}\circ M\in BMO$ is $\varepsilon_2$-sliceable
in the space $BMO$ and the martingale $\g\circ M\in BMO$ is
$\varepsilon_3$-sliceable in the space $BMO$.  Set
$\overline{C}_p:= q(1+C_p)+C_p$.  Let \be \label{contraction2}
\rho_2:=\overline{C}_p \max\left\{ \sqrt{2}p\, \varepsilon_3,\ \
2p\left\|\a\circ N_1\right\|_{BMO}\varepsilon_1+2p\,
\varepsilon_2^2\right\}<1.\ee

Then for any $(\xi,J)\in (L^p(\cF_T))^n\times (\cR^p)^n$, the
BSDE~(\ref{nbsde}) has a unique solution $(Y,Z\circ M, M^\bot)\in
(\cR^p)^n\times (\cH^p)^{2n}$. Moreover, there is a universal
constant $K_p$, which is independent of $(\xi,J)$,  such that
\be\label{bsde apriori} \left\|Y\right\|_{(\cR^p)^n}+\left\|(M,
M^\bot)\right\|_{(\cH^p)^{2n}}\le
K_p\left[\|\xi\|_{(L^p)^n}+\|J\|_{(\cR^p)^n}\right]. \ee\et

{\bf Proof of Theorem~\ref{bsde-basic}. } We shall still use the
contraction mapping principle and look for a fix-point. Consider
the following map $I$ in the Banach space $(\cR^p)^n\times
(\cH^p)^n$: for $(y,z\circ M)\in (\cR^p)^n\times (\cH^p)^n$,
define $I(y,z\circ M)$ to be components $(Y,Z\circ M)$ of the
unique adapted solution $(Y,Z\circ M,M^\bot)$ of the following
BSDE: \be\left\{\ba{rcl} Y_t&=&\displaystyle
\xi+J_T-J_t+\int_t^Tf(s, y_s)\, d\langle
N_1,N_2\rangle_s+\int_t^Tg(s,y_s,z_s)\, d\langle
M\rangle_s\\[0.3cm]
&&\displaystyle -\int_t^TZ_s\, dM_s-\int_t^T dM^\bot_s, \quad t\in
[0,T];\quad \langle M, M^\bot\rangle=0. \ea \right. \ee

We have \be \ba{rcl} Y_t&=& \displaystyle
E\left[\xi+(J_T-J_t)+\int_t^Tf(s, y_s)\, d\langle
N_1,N_2\rangle_s+\int_t^Tg(s,y_s,z_s)\, d\langle M\rangle_s\biggm
|\cF_t\right]\\[0.3cm]
 &=&\displaystyle -J_t+ E\left[\xi+J_T\
|\cF_t\right]+E\left[\int_t^Tf(s, y_s)\, d\langle
N_1,N_2\rangle_s\biggm|\cF_t\right]\\
&& \displaystyle +E\left[\int_t^Tg(s,y_s,z_s)\, d\langle
M\rangle_s\biggm |\cF_t\right].\ea\ee In view of  Doob's
inequality, we have \be \ba{rcl}
\|Y\|_{(\cR^p)^n}&\le&\displaystyle
\|J\|_{(\cR^p)^n}+\left\{E\left[\max_{0\le t\le
T}\left|E\left[\xi+J_T\, |\cF_t\right]\right|^p\, \right]\right\}^{1/p}\\[0.3cm]
&&\displaystyle+\left\{E\left[\max_{0\le t\le
T}\left|E\left[\int_t^Tf(s, y_s)\, d\langle
N_1,N_2\rangle_s\biggm|\cF_t\right]\right|^p\, \right]\right\}^{1/p}\\[0.4cm]
&&\displaystyle+\left\{E\left[\max_{0\le t\le
T}\left|E\left[\int_t^Tg(s,y_s,z_s)\,
d\langle M\rangle_s\biggm |\cF_t\right]\right|^p\, \right]\right\}^{1/p}\\[0.5cm]
&\le&\displaystyle
\|J\|_{(\cR^p)^n}+q\left\|\xi+J_T\right\|_{(L^p)^n}\\[0.3cm]
&&\displaystyle+\left\{E\left[\max_{0\le t\le
T}\left(E\left[\int_0^T|f(s, y_s)|\, |d\langle
N_1,N_2\rangle_s|\biggm|\cF_t\right]\right)^p\, \right]\right\}^{1/p}\\[0.4cm]
&&\displaystyle+\left\{E\left[\max_{0\le t\le
T}\left(E\left[\int_0^T|g(s,y_s,z_s)|\,
d\langle M\rangle_s\biggm |\cF_t\right]\right)^p\, \right]\right\}^{1/p}\\[0.4cm]
&\le&\displaystyle
\|J\|_{(\cR^p)^n}+q\left\|\xi+J_T\right\|_{(L^p)^n}+q\left\{E\left[\left(\int_0^T|f(s,
y_s)|\, |d\langle
N_1,N_2\rangle_s|\right)^p\, \right]\right\}^{1/p}\\[0.4cm]
&&\displaystyle +q\left\{E\left[\left(\int_0^T|g(s,y_s,z_s)|\,
d\langle M\rangle_s\right)^p\, \right]\right\}^{1/p}. \ea \ee
Proceeding identically as in the derivation of inequality
(\ref{Est1})  in the proof of Theorem~\ref{nse}, we have \be
E\left[\left|\int_0^T|f(s, y_s)|\, |d\langle
N_1,N_2\rangle_s|\right|^p\, \right] \le (\sqrt{2}p)^p\|\a\circ
N_1\|^p_{BMO}\|y\circ N_2\|^p_{(\cH^p)^n}.\ee Proceeding similarly
as in the derivation of inequality (\ref{Est2}), using the
Lipschitz assumption on $g$, we have
 \be\ba{rcl}&&
\displaystyle \left\{E\left[\left|\int_0^T|g(s,y_s,z_s)|\,
d\langle M\rangle_s\right|^p\, \right]\right\}^{1/p}\\
&\le & \displaystyle\left \{E\left[\left|\int_0^T(\b_s|y_s|+\g_s
|z_s|)\,
d\langle M\rangle_s\right|^p\, \right]\right\}^{1/p} \\[0.3cm]
&= & \displaystyle \left\|\langle \sqrt{\b}\circ M,
\sqrt{\b}|y|\circ M\rangle_T+\langle \g\circ M, |z|\circ
M\rangle_T\right\|_{L^p}.
 \ea \ee
Therefore, we have
  \be \ba{rcl} \|Y\|_{(\cR^p)^n}
&\le & \displaystyle
\|J\|_{(\cR^p)^n}+q\left\|\xi+J_T\right\|_{(L^p)^n}
+\sqrt{2}pq\|\a\circ
N_1\|_{BMO}\|y\circ N_2\|_{(\cH^p)^n} \\[0.3cm]
&& \displaystyle+q\left\|\langle \sqrt{\b}\circ M,
\sqrt{\b}|y|\circ M\rangle_T+\langle \g\circ M, |z|\circ
M\rangle_T\right\|_{L^p}\\[0.3cm]
&\le & \displaystyle
\|J\|_{(\cR^p)^n}+q\left\|\xi+J_T\right\|_{(L^p)^n}+\sqrt{2}pq\left\|\a\circ
N_1\|_{BMO}\|y\circ
N_2\right\|_{(\cH^p)^n}\\[0.3cm]
&& \displaystyle+\sqrt{2}pq\left \|\sqrt{\b}\circ
M\right\|_{BMO}\left\|\sqrt{\b}|y|\circ
M\right\|_{\cH^p}\\[0.3cm]
&&\displaystyle +\sqrt{2}pq\left\|\g\circ
M\right\|_{BMO}\left\|z\circ M\right\|_{(\cH^p)^n}\\[0.3cm]
&&\qquad \hbox{\rm (using Lemma~\ref{Lp}) } \\[0.3cm]
&\le &\displaystyle
\|J\|_{(\cR^p)^n}+q\left\|\xi+J_T\right\|_{(L^p)^n}+2pq\left\|\a\circ
N_1\right\|_{BMO}\left\|N_2\right\|_{BMO}\|y\|_{(\cR^p)^n}\\[0.3cm]
&& \displaystyle+2 pq\left\|\sqrt{\b}\circ
M\right\|_{BMO}^2\|y\|_{(\cR^p)^n} +\sqrt{2}pq\left\|\g\circ
M\right\|_{BMO}\left\|z\circ
M\right\|_{(\cH^p)^n} \\[0.3cm]
&&\qquad \qquad \quad \hbox{\rm (using Lemma~\ref{Emery})
}\\[0.3cm]
&\le&\displaystyle
\|J\|_{(\cR^p)^n}+q\left\|\xi+J_T\right\|_{(L^p)^n}+\sqrt{2}pq\left\|\g\circ
M\right\|_{BMO}\left\|z\circ
M\right\|_{(\cH^p)^n}\\[0.3cm]
&& \displaystyle+2pq\left(\left\|\a\circ
N_1\right\|_{BMO}\left\|N_2\right\|_{BMO}+\left\|\sqrt{\b}\circ
M\right\|_{BMO}^2\right)\|y\|_{(\cR^p)^n}. \ea\ee

Further, we have \be \ba{rcl}&&\displaystyle  \int_t^TZ_s\,
dM_s+\int_t^T dM^\bot_s\\
&=&\displaystyle \xi+J_T-J_t-Y_t +\int_t^Tf(s, y_s)\, d\langle
N_1,N_2\rangle_s\\
&=&\displaystyle+\int_t^Tg(s,y_s,z_s)\, d\langle M\rangle_s, \quad
t\in [0,T];\quad \langle M, M^\bot\rangle=0. \ea\ee From the BDG
inequality and using the similar arguments to the above, we have
\be \ba{rcl} &&\displaystyle \|z\circ M\|_{(\cH^p)^n}\le
C_p\left\|z\circ M+M^\bot\right\|_{(\cR^p)^n}\\[0.3cm]
&\le & \displaystyle
C_p\|\xi+J_T\|_{(L^p)^n}+C_p\|J\|_{(\cR^p)^n}+C_p\|Y\|_{(\cR^p)^n}\\[0.3cm]
&& \displaystyle+C_p\left\{E\left[\max_{0\le t\le
T}\left|\int_t^Tf(s, y_s)\, d\langle
N_1,N_2\rangle_s\right|^p\, \right]\right\}^{1/p}\\[0.3cm]
&&\displaystyle +C_p\left\{E\left[\max_{0\le t\le
T}\left|\int_t^Tg(s,y_s,z_s)\,
d\langle M\rangle_s\right|^p\, \right]\right\}^{1/p}\\[0.5cm]
&\le & \displaystyle
C_p\|\xi+J_T\|_{(L^p)^n}+C_p\|J\|_{(\cR^p)^n}+C_p\|Y\|_{(\cR^p)^n}\\[0.3cm]
&& \displaystyle+C_p\left\{E\left[\left(\int_0^T|f(s, y_s)|\,
|d\langle
N_1,N_2\rangle_s|\right)^p\, \right]\right\}^{1/p}\\[0.3cm]
&&\displaystyle +C_p\left\{E\left[\left(\int_0^T|g(s,y_s,z_s)|\,
d\langle M\rangle_s\right)^p\, \right]\right\}^{1/p}\\[0.5cm]
&\le & \displaystyle
C_p\|\xi+J_T\|_{(L^p)^n}+C_p\|J\|_{(\cR^p)^n}+C_p\|Y\|_{(\cR^p)^n}\\[0.3cm]
&& \displaystyle+\sqrt{2}pC_p\left\|\a\circ
N_1\right\|_{BMO}\left\|y\circ N_2\right\|_{(\cH^p)^n}\\[0.3cm]
&&\displaystyle +\sqrt{2}pC_p\left\|\sqrt{\b}\circ
M\right\|_{BMO}\left\|\sqrt{\b}|y|\circ
M\right\|_{\cH^p}\\[0.3cm]
& & \displaystyle+\sqrt{2}pC_p\left\|\g\circ
M\right\|_{BMO}\left\|\,|z|\circ M\, \right\|_{\cH^p}\\[0.3cm]
&\le & \displaystyle C_p\|\xi+J_T\|_{(L^p)^n}+C_p\|J\|_{(\cR^p)^n}+C_p\|Y\|_{(\cR^p)^n}\\[0.3cm]
&& \displaystyle+2pC_p\left\|\a\circ
N_1\right\|_{BMO}\left\|N_2\right\|_{BMO}\|y\|_{(\cR^p)^n}\\[0.3cm]
 &&\displaystyle
+2pC_p\left\|\sqrt{\b}\circ
M\right\|^2_{BMO}\|y\|_{(\cR^p)^n}+\sqrt{2}pC_p\left\|\g\circ
M\right\|_{BMO}\left\|\,|z|\circ M\,\right\|_{\cH^p}.\ea \ee

Concluding the above, we have \be \ba{rcl} &&\displaystyle
 \|Y\|_{(\cR^p)^n}+\|z\circ M\|_{(\cH^p)^n} \\
 &\le & \displaystyle (1+C_p)\|Y\|_{(\cR^p)^n}+C_p\|\xi+J_T\|_{(L^p)^n}+C_p\|J\|_{(\cR^p)^n}\\[0.3cm]
&& \displaystyle+2pC_p\left\|\a\circ
N_1\right\|_{BMO}\left\|N_2\right\|_{BMO}\|y\|_{(\cR^p)^n}\\[0.3cm]
 &&\displaystyle
+2pC_p\left\|\sqrt{\b}\circ
M\right\|^2_{BMO}\|y\|_{(\cR^p)^n}+\sqrt{2}pC_p\left\|\g\circ
M\right\|_{BMO}\left\|\,|z|\circ M\,\right\|_{\cH^p}\\[0.3cm]
&\le & \displaystyle \overline{C}_p\|\xi+J_T\|_{(L^p)^n}+\left(1+2C_p\right)\|J\|_{(\cR^p)^n}\\[0.3cm]
&& \displaystyle+2p\overline{C}_p\left\|\a\circ
N_1\right\|_{BMO}\left\|N_2\right\|_{BMO}\|y\|_{(\cR^p)^n}\\[0.3cm]
 &&\displaystyle
+2p\overline{C}_p\left\|\sqrt{\b}\circ
M\right\|^2_{BMO}\|y\|_{(\cR^p)^n}+\sqrt{2}p\overline{C}_p\left\|\g\circ
M\right\|_{BMO}\left\|\,|z|\circ M\,\right\|_{\cH^p}.\ea\ee

Let $(y^i,z^i\circ M)\in (\cR^p)^n\times (\cH^p)^n$ with $i=1,2$.
Denote by $(Y^i,Z^i\circ M)$ the image $I(y^i,z^i\circ M)$ for
$i=1,2$. Similar to the above arguments, we can show that \be
\label{estimate2}\ba{rcl} &&\displaystyle
\left\|Y^1-Y^2\right\|_{(\cR^p)^n}+\left\|(Z^1-Z^2)\circ M\right\|_{(\cH^p)^n}\\[0.3cm]
&\le &\displaystyle 2p\overline{C}_p\left\|\a\circ
N_1\right\|_{BMO}\left\|y^1-y^2\right\|_{(\cR^p)^n} \left\|N_2\right\|_{BMO}\\[0.3cm]
&& \displaystyle + 2p\overline{C}_p\left\|\sqrt{\b}\circ
M\right\|^2_{BMO}\left\|y^1-y^2\right\|_{(\cR^p)^n}\\[0.3cm]
&&\displaystyle +\sqrt{2}p\overline{C}_p\left\|\g\circ
M\right\|_{BMO}\left\|(z^1-z^2)\circ
M\right\|_{(\cH^p)^n} \\[0.3cm]
&=& \displaystyle 2p\overline{C}_p\left[\left\|\a\circ
N_1\right\|_{BMO}\left\|N_2\right\|_{BMO}+\left\|\sqrt{\b}\circ
M\right\|^2_{BMO}\right]\left\|y^1-y^2\right\|_{(\cR^p)^n} \\[0.3cm]
&&  \displaystyle + \sqrt{2}p\overline{C}_p\left\|\g\circ
M\right\|_{BMO}\left\|(z^1-z^2)\circ M\right\|_{(\cH^p)^n}\\[0.3cm] &\le&
\max\left\{ \sqrt{2}p\left\|\g\circ M\right\|_{BMO}, \  2p
\left\|\a\circ N_1\right\|_{BMO}\left\|N_2\right\|_{BMO}+2p\left\|\sqrt{\b}\circ M\right\|^2_{BMO}\right\}\\[0.3cm]
&& \displaystyle \times \overline{C}_p
\left[\left\|y^1-y^2\right\|_{(\cR^p)^n}+\left\|(z^1-z^2)\circ
M\right\|_{(\cH^p)^n}\right].\ea\ee

Since the martingales $N_2\in BMO, \sqrt{\b}\circ M\in BMO$, and
$\g\circ M\in BMO$ are respectively $\varepsilon_1$-sliceable,
$\varepsilon_2$-sliceable, and $\varepsilon_3$-sliceable,
 there is a finite sequence of stopping times $\{T_i,
i=1,2,\cdots,{\widetilde I} \}$ such that the following are
satisfied:

(i) $0=T_0\le T_1 \le T_2\le \cdots\le T_{\widetilde I}\le
T_{\widetilde I+1}=T$;

(ii) $\left\|N_{2i}\right\|_{BMO}\le \varepsilon_1$,
$\left\|\sqrt{\b}\circ M_i\right\|_{BMO}\le \varepsilon_2$ and
$\left\|\b\circ M_i\right\|_{BMO}\le \varepsilon_3$ where
$N_{2i}:=N_2^{T_{i+1}}-N_2^{T_i}$ and $M_i:=M^{T_{i+1}}-M^{T_i}$
are defined on $[T_i,T_{i+1}]$.

Since $\left\|\a\circ N_{1i}\right\|_{BMO}\le \left\|\a\circ
N_1\right\|_{BMO}$, we have \be
\label{contraction22}\rho_{2i}:=\overline{C}_p\max\left\{
\sqrt{2}p\varepsilon_3,\ \  2p\left\|\a\circ
N_{1i}\right\|_{BMO}\varepsilon_1+2p\varepsilon_2^2\right\}\le
\rho_2.\ee with $N_{1i}:=N_1^{T_{i+1}}-N_1^{T_i}$ for $i=0,1,
2,\cdots,{\widetilde I}$.

Set $\cR^p_i:=\cR^p[T_i,T_{i+1}]$ and
$\cH^p_i:=\cH^p(T_i,T_{i+1})$ for $i=0,1,\cdots,\widetilde I$. Set
$Y^{\widetilde I+1}(T)=\xi$. Consider the map $I_i$ in the Banach
space $(\cR^p_i)^n\times (\cH^p_i)^n$: for $(y,z\circ M)\in
(\cR^p_i)^n\times (\cH^p_i)^n$, define $I_i(y,z\circ M)$ to be
components $(Y,Z\circ M)$ of the unique adapted solution
$(Y,Z\circ M,M^\bot)$ of the following BSDE: \be\left\{\ba{rcl}
Y_t&=&\displaystyle
Y^{i+1}_{T_{i+1}}+(J_{T_{i+1}}-J_t)+\int_t^{T_{i+1}}f(s, y_s)\,
d\langle N_{1i},N_{2i}\rangle_s+\int_t^{T_{i+1}}g(s,y_s,z_s)\,
d\langle
M_i\rangle_s\\[0.3cm]
&&\displaystyle -\int_t^{T_{i+1}}Z_s\, dM_{is}-\int_t^{T_{i+1}}
dM^\bot_s, \quad t\in [T_i,T_{i+1}];\quad \langle M,
M^\bot\rangle=0. \ea \right. \ee

Similar to the derivation of inequality~(\ref{estimate2}), we have
\be \ba{rcl} &&\displaystyle
\left\|I_i(y^1,z^1)-I_i(y^2,z^2)\right\|_{(\cR^p_i)^n\times (\cH^p_i)^n}\\[0.3cm]
 &\le&\displaystyle
\max\left\{ \sqrt{2}p\left\|\g\circ M_i\right\|_{BMO}, \  2p
\left\|\a\circ N_{1i}\right\|_{BMO}\left\|N_{2i}\right\|_{BMO}+2p\left\|\sqrt{\b}\circ M_i\right\|^2_{BMO}\right\}\\[0.3cm]
&& \displaystyle \times \overline{C}_p
\left[\left\|y^1-y^2\right\|_{(\cR^p_i)^n}+\left\|(z^1-z^2)\circ
M\right\|_{(\cH^p_i)^n}\right]\\[0.3cm]
&\le& \displaystyle \rho_{2i}
\left[\left\|y^1-y^2\right\|_{(\cR^p_i)^n}+\left\|(z^1-z^2)\circ
M\right\|_{(\cH^p_i)^n}\right]\\[0.3cm]
& \le & \displaystyle \rho_2
\left[\left\|y^1-y^2\right\|_{(\cR^p_i)^n}+\left\|(z^1-z^2)\circ
M\right\|_{(\cH^p_i)^n}\right]\ea\ee for any
$(y^1,z^1),(y^2,z^2)\in (\cR^p_i)^n\times (\cH^p_i)^n$.  In view
of inequality~(\ref{contraction2}) in the second assumption of the
theorem, we see that for each $i=0,1,\cdots,\widetilde I$, $I_i$
is a contraction map on $(\cR^p_i)^n\times (\cH^p_i)^n$. More
precisely, first, since $I_{\widetilde I}$ is a contraction, the
following BSDE: \be\left\{\ba{rcl} Y_t&=&\displaystyle
\xi+(J_T-J_t)+\int_t^{T}f(s, Y_s)\, d\langle N_{1\widetilde
I},N_{2\widetilde I}\rangle_s+\int_t^{T}g(s,Y_s,Z_s)\, d\langle
M_{\widetilde I}\rangle_s\\[0.3cm]
&&\displaystyle -\int_t^{T}Z_s\, dM_{{\widetilde I}s}-\int_t^{T}
dM^\bot_s, \quad t\in [T_{\widetilde I},T]\quad \langle M,
M^\bot\rangle=0 \ea \right. \ee has a unique solution
$(Y^{\widetilde I}, Z^{\widetilde I}\circ M_{\widetilde I},
M^{{\widetilde I}\bot})\in (\cR^p_{\widetilde I})^n\times
(\cH^p_{\widetilde I})^{2n}$. Second, consider the following BSDE:
\be\left\{\ba{rcl} Y_t&=&\displaystyle Y^{\widetilde
I}_{T_{\widetilde I}}+(J_{T_{\widetilde
I}}-J_t)+\int_t^{T_{\widetilde I}}f(s, Y_s)\, d\langle
N_{1,{\widetilde I}-1},N_{2, {\widetilde
I}-1}\rangle_s+\int_t^{T_{\widetilde I}}g(s,Y_s,Z_s)\, d\langle
M_{{\widetilde I}-1}\rangle_s\\[0.3cm]
&&\displaystyle -\int_t^{T_{\widetilde I}}Z_s\, dM_{{{\widetilde
I}-1},s}-\int_t^{T_{\widetilde I}} dM^\bot_s, \quad t\in
[T_{{\widetilde I}-1},T_{\widetilde I}];\quad \langle M,
M^\bot\rangle=0. \ea \right. \ee Since the map $I_{\widetilde
I-1}$ is a contraction in $(\cR^p_{\widetilde I})^n\times
(\cH^p_{\widetilde I})^{n}$, it has a unique solution
$(Y^{{\widetilde I}-1},Z^{{\widetilde I }-1}\circ M_{{\widetilde
I}-1},M^{{{\widetilde I}-1} \bot})$ in $(\cR^p_{{\widetilde
I}-1})^n\times (\cH^p_{{\widetilde I}-1})^{2n}$. Inductively in a
backward way, we can show that the following BSDE:
\be\left\{\ba{rcl} Y_t&=&\displaystyle
Y^{i+1}_{T_{i+1}}+(J_{T_{i+1}}-J_t)+\int_t^{T_{i+1}}f(s, Y_s)\,
d\langle N_{1i},N_{2i}\rangle_s+\int_t^{T_{i+1}}g(s,Y_s,Z_s)\,
d\langle
M_i\rangle_s\\[0.3cm]
&&\displaystyle -\int_t^{T_{i+1}}Z_s\, dM_{is}-\int_t^{T_{i+1}}
dM^\bot_s, \quad t\in [T_i,T_{i+1}];\quad \langle M,
M^\bot\rangle=0 \ea \right. \ee has a unique solution
$(Y^i,Z^i,M^{i,\bot})$ in $(\cR^p_i)^n\times (\cH_i^p)^{2n}$ for
$i=0,1,\cdots, \widetilde I$. Moreover, we have
\be\label{Est4}\ba{rcl} &&\displaystyle
\left(1-\rho_2\right)\left(\left\|Y^i\right\|_{(\cR^p_i)^n}+\left\|Z^i\right\|_{(\cH^p_i)^n}\right)\\[0.3cm]
&\le&\displaystyle \overline{C}_p\left
\|Y^{i+1}_{T_{i+1}}+J_{T_{i+1}}\right\|_{(L^p(\cF_{T_{i+1}}))^n}+(2C_p+1)\left\|J\right\|_{(\cR^p_i)^n}
\ea\ee for $i=0,1,\cdots, \widetilde I$.

Then, the triple of processes $(Y,Z\circ M,M^\bot)$ given by
\be\ba{rcl} X(t)&:=&\displaystyle \sum_{i=0}^{\widetilde
I}Y^i_t\chi_{[T_i,
T_{i+1})}(t), \quad t\in [0,T],\\
Z(t)&:=&\displaystyle \sum_{i=0}^{\widetilde I}Z^i_t\chi_{[T_i,
T_{i+1})}(t), \quad t\in [0,T],\\
M^\bot(t)&:=&\displaystyle
M^{0,\bot}_t\chi_{[0,T_1)}(t)+\sum_{i=1}^{\widetilde
I}\left[M^{i\bot}_t+M^{i-1,\bot}_{T_i}\right]\chi_{[T_i,
T_{i+1})}(t), \quad t\in [0,T]\\
\ea \ee lies in $(\cR^p)^n\times (\cH^p)^{2n}$ and is the unique
adapted  solution to BSDE~(\ref{nbsde}). The estimate~(\ref{bsde
apriori}) is a consequence of the inequalities~(\ref{Est4}).
\endpf

Consider BSDE~(\ref{nbsde}) for the case of $f=0, J=0$ and $g$
being independent of $y$. That is, consider the following
nonlinear BSDEs: \be\left\{\ba{rcl} Y_t&=&\displaystyle
\xi+\int_t^Tg(s,Z_s)\, d\langle
M\rangle_s\\[0.3cm]
&&\displaystyle -\int_t^TZ_s\, dM_s-\int_t^T dM^\bot_s, \quad t\in
[0,T];\quad \langle M, M^\bot\rangle=0. \ea
\right.\label{nbsdeinfty}\ee For the extremal case of $p=\infty$,
we have the following result.

\bt\label{bsde-basic-infty} Assume that

(i) There is an $\{\cF_t, 0\le t\le T\}$-adapted processes
 $\g(\cdot)$ such that  \be g(\cdot,0)=0; \quad
|g(t,z_1)-g(t,z_2)|\le  \g (t)|z_1-z_2|\ee for $z_1,z_2\in R^n.$

(ii) The martingale $\g\circ M\in BMO$ is $\varepsilon$-sliceable
in the space $BMO$  such that \be \label{contraction2infty}
\sqrt{2}\varepsilon<1.\ee

Then for $\xi\in (BMO)^n$, the BSDE~(\ref{nbsdeinfty}) has a
unique solution $(Y,Z\circ M, M^\bot)$ such that $ (Z\circ M,
M^\bot)\in (BMO)^{2n}$. Moreover, there is a universal constant
$K$, which is independent of $\xi$, such that \be \left\|(Z\circ
M, M^\bot)\right\|_{(BMO)^{2n}}\le K \|\xi\|_{(BMO)^n}.\ee \et

{\bf Proof of Theorem~\ref{bsde-basic-infty}. }
 We shall still use the
contraction mapping principle and look for a fix-point. Consider
the following map $I$ in the Banach space $(BMO)^n$: for $z\circ
M\in (BMO)^n$, define $I(z\circ M)$ to be component $Z\circ M$ of
the unique adapted solution $(Y,Z\circ M,M^\bot)$ of the following
BSDE: \be\left\{\ba{rcl} Y_t&=&\displaystyle
\xi+\int_t^Tg(s,z_s)\, d\langle
M\rangle_s\\[0.3cm]
&&\displaystyle -\int_t^TZ_s\, dM_s-\int_t^T dM^\bot_s, \quad t\in
[0,T];\quad \langle M, M^\bot\rangle=0. \ea \right. \ee The
following shows that $I(z\circ M)$ is in the BMO space for any
$z\circ M\in (BMO)^n$: \be\ba{rcl}&& \|Z\circ M\|_{(BMO)^n}\le
\left\|Z\circ
M+M^\bot\right\|_{(BMO)^n}\\[0.3cm]
&=& \displaystyle \left\|\xi+\int_0^Tg(s,z_s)\, d\langle M\rangle_s\right\|_{(BMO)^n}\\
&\le&\displaystyle  \|\xi\|_{(BMO)^n}+\left\|\int_0^Tg(s,z_s)\,
d\langle M\rangle_s\right\|_{(BMO)^n}\\[0.3cm]
&\le&\displaystyle  \|\xi\|_{(BMO)^n}+\left\|\int_0^T\g_s z_s\,
d\langle M\rangle_s\right\|_{(BMO)^n}\\[0.3cm]
&=&\displaystyle  \|\xi\|_{(BMO)^n}+\left\|\langle \g\circ M,
z\circ
M\rangle_T\right\|_{(BMO)^n}\\[0.3cm]
&\le&\displaystyle  \|\xi\|_{(BMO)^n}+ \sqrt{2}\|\g\circ M\|_{BMO}
\|z\circ M\|_{(BMO)^n}\\[0.3cm]
&& \displaystyle \qquad \qquad\qquad \hbox{\rm (using
Lemma~\ref{Linfty})}.
 \ea\ee

Let $z^i\circ M\in (BMO)^n$ with $i=1,2$. Denote by $Z^i\circ M$
the image $I(z^i\circ M)$ for $i=1,2$. Similar to the above
arguments, we can show that \be \ba{rcl} && \left\|Z^1\circ
M-Z^2\circ M\right\|_{(BMO)^n}\le \sqrt{2}\left\|\g\circ
M\right\|_{BMO} \left\|z^1\circ M-z^2\circ M\right\|_{(BMO)^n}.
\ea\ee

The rest of the proof is identical to that of
Theorem~\ref{bsde-basic}.
\endpf

\br Theorem~\ref{bsde-basic-infty} is not implied by
Theorem~\ref{bsde-basic} due to the fact that the assumption (ii)
of the latter involves $p$. In fact, the proof of the former
appeals to Lemma~\ref{Linfty}, while the proof of the latter
appeals to  Lemma~\ref{Lp}. Lemma~\ref{Linfty} is not implied by
Lemma~\ref{Lp}.\er

\subsection{ Comments on the slice-ability assumption in the space $BMO$ on the
martingales $N_2, \g\circ M$, and $\b\circ M$ in
Theorems~\ref{nse}, ~\ref{bsde-basic}, and~\ref{bsde-basic-infty}}

 Schachermayer~\cite{Schacher} shows that any martingale in $\overline {\cH^\infty}^{BMO}$ is sliceable in the
space $BMO$. Therefore, the suitable slice-ability assumption in
the space $BMO$ in the preceding subsection on the martingales
$N_2, \g\circ M$, and $\b\circ M\in BMO$ is automatically true
when they are in the space $\overline {\cH^\infty}^{BMO}$.
Therefore, we have the following

\bt \label{nsesliceable} Let $p\in [1, \infty)$. Assume that

(i) There are two $\{\cF_t, 0\le t\le T\}$-adapted processes
$\a(\cdot)$ and $\b(\cdot)$ such that \be f(t,0)=0; \quad
|f(t,x_1)-f(t,x_2)|\le \a (t)\max_{0\le s\le t}|x_1(s)-x_2(s)|,
\quad x_1,x_2\in \cD^n\ee and \be g(t,0)=0; \quad
|g(t,x_1)-g(t,x_2)|\le \b (t)\max_{0\le s\le t}|x_1(s)-x_2(s)|,
\quad x_1,x_2\in \cD^n.\ee

(ii) The martingale  $\a \circ N_1\in BMO$. Both martingales $
N_2$ and $\b\circ M$ are in the space
${\overline{\cH^\infty}}^{BMO}$.

Then for any $J\in (\cR^p)^n$, there is unique solution in
$(\cR^p)^n$ to equation~(\ref{nsdeG}). Furthermore, there is a
constant $K_p$, which is independent of $J$, such that \be
\label{aprioriinqualsliceable} \|X\|_{\cR^p}\le K_p
\|J\|_{\cR^p}.\ee If $J\in (\cR^p)^n$ is a semi-martingale, then
so  is the solution. \et

\br Theorem~\ref{nsesliceable} more or less generalizes
Protter~\cite[Lemma 2, page 252]{Protter1}.\er

\bc\label{cor} There are three real valued nonnegative $\{\cF_t,
0\le t\le T\}$-adapted processes $\a(\cdot), \b(\cdot)$ and
$\g(\cdot)$ such that the adapted $R^{n\times n}$-valued processes
$A,B$, and $D$ are bounded respectively by $\a,\b$, and $\g$.
Assume that the martingale  $\a \circ N_1\in BMO$, and the
martingales $ N_2, \sqrt{\b}\circ M$, and $\g\circ M$ are all in
the space ${\overline{\cH^\infty}}^{BMO}$. Let $S(\cdot)$ be the
fundamental solution matrix process to the following SDE:
\be\label{hsde}\left\{ \ba{rcl}   dS(t)&=&\displaystyle \left[
A_t\, d\langle N_1, N_2\rangle_t+B_t\,
d\langle M\rangle_t +D_t\, dM_t\right]S(t), \quad t\in [0,T];\\
S(0)&=&I.\ea\right.\ee Then, for any $p\in [1, \infty)$, there is
a universal constant $K_p$ such that for any stopping time $\tau$,
we have  \be\label{univ} E\left[\max_{\tau\le t\le
T}\left|S^{-1}(\tau)S(t)\right|^p\,\biggm |\cF_\tau \right]\le
K_p^p. \ee The last inequality implies that $S(\cdot)$ satisfies
the reverse H\"older property $(R_p)$ for any $p\in [1,
\infty)$.\ec

{\bf Proof of Corollary~\ref{cor}. } The assumptions of
Theorem~\ref{nsesliceable} are all satisfied except that the two
continuous local martingales $N_1,N_2$, and the real nonnegative
process $\a$ in Theorem~\ref{nsesliceable} correspond to two
two-dimensional vector-valued continuous local  martingales $(N_1,
M), (N_2, \sqrt{\b}\circ M)$, and the two-dimensional
vector-valued
 processes $(\a, \sqrt{\b})$ in this corollary.

 Consider any stopping time $\tau$. Take any
$G\in \cF_\tau.$ For $J=S(\tau)\chi_{G}$, it is easy to see that
$X:=S\chi_{G}$ is the unique solution to the SDE~(\ref{hsde}) with
the initial condition being replaced with
$X(\tau)=S(\tau)\chi_{G}$. The assertions of
Theorem~\ref{nsesliceable} are still true for $X$. In view of the
estimate~(\ref{aprioriinqualsliceable}) of
Theorem~\ref{nsesliceable}, we have \be E\left[\max_{\tau\le t\le
T}\left|S(t)\chi_{G}\right|^p\, \right]\le
K_p^pE\left[\left|S(\tau)\chi_{G}\right|^p\, \right].\ee
Therefore, we have \be E\left[\max_{\tau\le t\le
T}\left|S(t)\right|^p\chi_{G}\right]\le K_p^p
E\left[\left|S(\tau)\right|^p\chi_{G}\right].\ee This implies the
inequality~(\ref{univ}).
 The proof is complete.
\endpf

\bt\label{bsde-basic-sliceable} Let $p\in (1, \infty)$.  Assume
that

(i) There are three $\{\cF_t, 0\le t\le T\}$-adapted processes
$\a(\cdot), \b(\cdot)$ and $\g(\cdot)$ such that \be f(\cdot,0)=0;
\quad |f(t,y_1)-f(t,y_2)|\le \a (t)|y_1-y_2|\ee for $ y_1,y_2\in
R^n$ and \be g(\cdot,0,0)=0; \quad |g(t,y_1,z_1)-g(t,y_1,z_2)|\le
\b (t)|y_1-y_2|+ \g (t)|z_1-z_2|\ee for $y_1,y_2,z_1,z_2\in R^n.$

(ii) The martingale  $\a \circ N_1\in BMO$. The martingales $ N_2,
\sqrt{\b}\circ M$, and $\g\circ M$ are all in the space
${\overline{\cH^\infty}}^{BMO}$.

Then for $(\xi, J)\in (L^p(\cF_T))^n\times (\cR^p)^n$,
BSDE~(\ref{nbsde}) has a unique solution $(Y,Z\circ M, M^\bot)\in
(\cR^p)^n\times (\cH^p)^{2n}$. Moreover, there is a universal
constant $K_p$, which is independent of $J$,  such that \be
\left\|Y\right\|_{(\cR^p)^n}+\left\|(M,
M^\bot)\right\|_{(\cH^p)^{2n}}\le
K_p\left[\|\xi\|_{(L^p)^n}+\|J\|_{(\cR^p)^n}\right]. \ee\et

Note that the existence and uniqueness of F\"ollmer-Schweizer
decomposition (see F\"ollmer and Schweizer~\cite{FSch}) is exactly
the existence and uniqueness of a one-dimensional linear BSDE, but
possibly and typically with unbounded coefficients.
Theorem~\ref{bsde-basic-sliceable} includes as particular cases
the existence and uniqueness results on linear BSDEs not only for
bounded coefficients by Bismut~\cite{Bis3}, but also for unbounded
coefficients  by Monat and
Stricker~\cite{MonatStricker1,MonatStricker2,MonatStricker3} and
by Schweizer~\cite{Schweizer1,Schweizer2}---where $\gamma\circ M$
is assumed to be in $\cH^\infty$---in the case of no jumps in
$\gamma\circ M$, and by Delbaen et al~\cite{DMSSS1,DMSSS2} in the
case of $\gamma\circ M \in {\overline{\cH^\infty}}^{BMO}$. Note
that $ \lambda\circ M \in {\overline{\cH^\infty}}^{BMO}$ when the
process $\lambda$ is a uniformly bounded adapted process and the
local martingale $M$ is a Brownian motion stopped at a finite
deterministic time $T$. Therefore,
Theorem~\ref{bsde-basic-sliceable} also includes as particular
cases the existence and uniqueness results on nonlinear BSDEs of
Pardoux and Peng~\cite{PP} (for $L^2$ solutions), and El Karoui,
Peng, and Quenez~\cite[Theorem 5.1, page 54]{ElPengQuenez} (for
$L^p$  solutions ($p>1$)).

\bt\label{bsde-basic-infty-sliceable} Assume that

(i) There is an $\{\cF_t, 0\le t\le T\}$-adapted processes
 $\g(\cdot)$ such that  \be \int_0^T|g(t,0)||d\langle M\rangle_s|\in BMO; \quad
|g(t,z_1)-g(t,z_2)|\le  \g (t)|z_1-z_2|\ee for $z_1,z_2\in R^n.$

(ii) The martingale $\g\circ M\in {\overline{\cH^\infty}}^{BMO}$.

Then for $\xi\in (BMO)^n$, BSDE~(\ref{nbsdeinfty}) has a unique
solution $(Y,Z\circ M, M^\bot)$ such that $ (Z\circ M, M^\bot)\in
(BMO)^{2n}$. Moreover, there is a universal constant $K$, which is
independent of $\xi$, such that \be \left\|(Z\circ M,
M^\bot)\right\|_{(BMO)^{2n}}\le K \|\xi\|_{(BMO)^n}.\ee \et

When the data (i.e., the terminal state $\xi$ and the ``zero" term
$J$) is essentially bounded, instead of just being in the BMO
space, the unique adapted solution $(Y, Z\circ M, M^\bot)$ to
BSDE~(\ref{nbsde}) can be further proved to lie in the better
space: $(\cR^\infty)^n\times
\left({\overline{\cH^\infty}}^{BMO}\right)^n$.

\bt \label{nbsdebounded} Assume that

(i)  There are three $\{\cF_t, 0\le t\le T\}$-adapted processes
$\a(\cdot), \b(\cdot)$ and $\g(\cdot)$ such that \be f(\cdot,0)=0;
\quad |f(t,y_1)-f(t,y_2)|\le \a (t)|y_1-y_2|\ee for $ y_1,y_2\in
R^n$ and \be g(\cdot,0,0)=0; \quad |g(t,y_1,z_1)-g(t,y_1,z_2)|\le
\b (t)|y_1-y_2|+ \g (t)|z_1-z_2|\ee for $y_1,y_2,z_1,z_2\in R^n.$

(ii) The martingale  $\a \circ N_1\in BMO$. The martingales $ N_2,
\sqrt{\b}\circ M$, and $\g\circ M$ are all in the space
${\overline{\cH^\infty}}^{BMO}$.

For any $\xi \in L^\infty(\cF_T)$ and $J
 \in \cR^\infty$,   there is
unique adapted solution $(Y, Z\circ M, M^\bot)$ to
BSDE~(\ref{nbsde}), with $Y\in (\cL^\infty)^n$ and $Z\circ M+
M^\bot\in \left(\overline{L^\infty}^{BMO}\right)^n$.  \et

{\bf Proof of Theorem~\ref{nbsdebounded}. } From
Theorem~\ref{bsde-basic-sliceable}, we know that there is a unique
adapted solution $(Y, Z\circ M, M^\bot)\in (\cR^p)^n\times
(\cH^p)^{2n}$ to BSDE~(\ref{nbsde}) for any $p\in (1,\infty)$. The
proof is divided into the following three steps.

Step 1. We show that $Y\in (\cL^\infty)^n$. In fact,
BSDE~(\ref{nbsde}) can be written into the following linear form:
\be \label{newform}\ba{rcl} Y_t&=&\displaystyle
\xi+J_T-J_t+\int_t^T A_s^\tau Y_s \, d\langle N_1, N_2\rangle_s
+\int_t^T\left(D_s^\tau Z_s+B_s^\tau Y_s\right)\, d\langle M\rangle_s \\
&& \displaystyle -\int_t^TZ_s\, dM_s-\int_t^T\, dM^\bot_s,  \quad
t\in [0,T]; \quad \langle M, M^\bot\rangle=0,\ea\ee with the
adapted matrix-valued processes $A,B$, and $D$ being bounded
respectively by $\a,\b$, and $\g$. Let $S(\cdot)$ be the
fundamental solution matrix process to the SDE~(\ref{hsde}). Then,
we have \be Y_t=E\left[S^\tau
(t)^{-1}S^\tau(T)(\xi+J_T)-\int_t^TS^\tau(t)^{-1}S^\tau(s)[A_s^\tau
J_s \, d\langle N_1, N_2\rangle_s+B_s^\tau J_s \, d\langle
M\rangle_s]\biggm |\cF_t\right]-J_t.\ee In view of
Corollary~\ref{cor}, $S(\cdot)$ satisfies the reverse H\"older
property $(R_p)$ for any $p\in [1, \infty)$, and the
inequality~(\ref{univ}) hold. Therefore, we have \be \ba{rcl}
|Y_t|&\le& \displaystyle |J_t|+E\left[|S^\tau
(t)^{-1}S^\tau(T)|\cdot|\xi+J_T|\,
|\cF_t\right]\\[0.3cm]
&&\displaystyle +\|J\|_{(\cR^\infty)^n}E\left[\left(\sup_{s\in
[t,T]}|S^\tau(t)^{-1}S^\tau(s)|\right)\int_t^T[|A_s| \, |d\langle
N_1, N_2\rangle_s|+|B_s| \, d\langle M\rangle_s]\biggm
|\cF_t\right]\\[0.3cm]
&\le& \displaystyle |J_t|+E\left[|S^\tau
(t)^{-1}S^\tau(T)|\cdot|\xi+J_T|\,
|\cF_t\right]\\[0.3cm]
&&\displaystyle +\|J\|_{(\cR^\infty)^n}E\left[\left(\sup_{s\in
[t,T]}|S^\tau(t)^{-1}S^\tau(s)|\right)\int_t^T[\a_s\, |d\langle
N_1, N_2\rangle_s|+\b_s \, d\langle M\rangle_s]\biggm
|\cF_t\right]\\[0.3cm]
&\le & \displaystyle \|J\|_{(\cR^\infty)^n}+K_1\left(\|\xi\|_{(L^\infty)^n}+\|J\|_{(\cR^\infty)^n}\right)\\[0.3cm]
&&\displaystyle +K_2\|J\|_{(\cR^\infty)^n}\left\{E\left[
\left|\int_t^T\left(\, |d\langle \a\circ N_1, N_2\rangle_s|+
d\langle \sqrt{\b}\circ M\rangle_s\right)\right|^2\biggm
|\cF_t\right]
\right\}^{1/2}\\[0.3cm]
&\le & \displaystyle \|J\|_{(\cR^\infty)^n}+K_1\left(\|\xi\|_{(L^\infty)^n}+\|J\|_{(\cR^\infty)^n}\right)\\[0.3cm]
&&\displaystyle +K_2\|J\|_{(\cR^\infty)^n}\left\{E\left[
2\left|\langle \a\circ N_1, N_2\rangle_t^T\right|^2+
2\left|\langle \sqrt{\b}\circ M\rangle_t^T\right|^2\biggm
|\cF_t\right] \right\}^{1/2} .\ea\ee Here, $K_1$ and $K_2$ are
introduced in Corollary~(\ref{cor}).  In view of the assumption
(ii) of the theorem, using Kazamaki~\cite[Lemma 2.6, page 48]{Kaz}
and the John-Nirenberg
 inequality (see Kazamaki~\cite[Theorem 2.2, page 29]{Kaz}),
 we have $b(N_2)=b(\sqrt{\b}\circ
M)=\infty$ and  \be E\left[ \exp{\left(\epsilon\langle \a\circ
N_1\rangle_t^T\right)}\biggm |\cF_t\right]\in \cL^\infty\ee for
$\epsilon < \|\a\circ N_1\|^{-2}_{(BMO)^n}$. Then, the following
process \be \ba{rcl} && \displaystyle E\left[ 2\left|\langle
\a\circ N_1, N_2\rangle_t^T\right|^2+ 2\left|\langle
\sqrt{\b}\circ
M\rangle_t^T\right|^2\biggm |\cF_t\right]\\[0.3cm]
&\le& \displaystyle E\left[ 2\langle \a\circ N_1\rangle_t^T
\langle N_2\rangle_t^T+ 2\left|\langle \sqrt{\b}\circ
M\rangle_t^T\right|^2\biggm |\cF_t\right]\\[0.3cm]
&\le& \displaystyle E\left[ \left|\epsilon\langle \a\circ
N_1\rangle_t^T\right|^2+\epsilon^{-2}\left| \langle
N_2\rangle_t^T\right|^2+ 2\left|\langle \sqrt{\b}\circ
M\rangle_t^T\right|^2\biggm |\cF_t\right]\\[0.3cm]
&\le& \displaystyle E\left[ \left|\epsilon\langle \a\circ
N_1\rangle_t^T\right|^2\biggm
|\cF_t\right]+\epsilon^{-2}E\left[\left| \langle
N_2\rangle_t^T\right|^2\biggm |\cF_t\right]+ 2E\left[\left|\langle
\sqrt{\b}\circ M\rangle_t^T\right|^2\biggm |\cF_t\right]\\[0.3cm]
&\le& \displaystyle E\left[ \exp{\left(\epsilon\langle \a\circ
N_1\rangle_t^T\right)}\biggm
|\cF_t\right]+\epsilon^{-2}E\left[\exp{\left( \langle
N_2\rangle_t^T\right)}\biggm |\cF_t\right]\\[0.3cm]
&& \displaystyle+ 2E\left[\exp{\left(\langle \sqrt{\b}\circ
M\rangle_t^T\right)}\biggm |\cF_t\right]\in \cL^\infty\ea\ee for
$\epsilon \in(0, \|\a\circ N_1\|^{-2}_{(BMO)^n})$. Consequently,
we have $Y\in (\cL^\infty)^n$.

Step 2. We  show that $Z\circ M+ M^\bot\in BMO$. To simplify the
exposition, set \be
C_{YJ}:=\|Y\|_{(\cR^\infty)^n}+\|J\|_{(\cR^\infty)^n}.\ee In view
of BSDE~(\ref{newform}), using It\^o's formula and standard
arguments, we can obtain the following estimate for any stopping
time $\sigma$: \be\ba{rcl} &&\displaystyle E\left[ \langle Z\circ
M+
M^\bot\rangle_\sigma^T|\cF_\sigma\right] \\[0.3cm]
&\le &\displaystyle
E\left[|\xi+J_T|^2+2\int_\sigma^T|Y_s+J_s|\left[\a_s |d\langle
N_1, N_2\rangle_s|+(\b_s|Y_s|+\g_s|Z_s|)\, d\langle
M\rangle_s\right]
\, \biggm |\cF_\sigma \right]\\[0.3cm]
&\le &\displaystyle
2\|\xi\|_{L^\infty}^2+2\|J\|_{\cR^\infty}^2+2C_{YJ}
E\left[\int_\sigma^T\a_s |d\langle N_1,
N_2\rangle_s|+\int_\sigma^T(\b_s|Y_s|+\g_s|Z_s|)\, d\langle
M\rangle_s\, \biggm |\cF_\sigma \right]\\[0.3cm]
&\le &\displaystyle
2\|\xi\|_{L^\infty}^2+2\|J\|_{\cR^\infty}^2+2C_{YJ}E\left[\int_\sigma^T|d\langle
\a \circ N_1, N_2\rangle_s|\,
\biggm |\cF_\sigma \right]\\[0.3cm]
&&\displaystyle +2C_{YJ}\|Y\|_{(\cR^\infty)^n}E\left[\langle
\sqrt{\b}\circ M\rangle^T_\sigma \, \biggm |\cF_\sigma \right]
+2C_{YJ}E\left[\langle \g\circ M, |Z|\circ M\rangle_\sigma^T\, \biggm |\cF_\sigma \right]\\[0.3cm]
&\le &\displaystyle
2\|\xi\|_{L^\infty}^2+2\|J\|_{\cR^\infty}^2+2C_{YJ}E\left[\left(\langle
\a \circ N_1\rangle_\sigma^T\right)^{1/2} \left(\langle
N_2\rangle_\sigma^T\right)^{1/2} |\,
\biggm |\cF_\sigma \right]\\[0.3cm]
&&\displaystyle +2C_{YJ}\|Y\|_{(\cR^\infty)^n}E\left[\langle
\sqrt{\b}\circ M\rangle^T_\sigma \, \biggm |\cF_\sigma \right]\\[0.3cm]
&&\displaystyle+2C_{YJ}E\left[\left(\langle \g\circ
M\rangle_\sigma^T\right)^{1/2}\left(\langle |Z|\circ
M\rangle_\sigma^T\right)^{1/2}\, \biggm |\cF_\sigma \right]\\[0.3cm]
&\le &\displaystyle
2\|\xi\|_{L^\infty}^2+2\|J\|_{\cR^\infty}^2+2C_{YJ}\|\a \circ
N_1\|_{BMO} \|
N_2\|_{BMO}\\[0.3cm]
&&\displaystyle +2C_{YJ}\|Y\|_{(\cR^\infty)^n}\|
\sqrt{\b}\circ M\|^2_{BMO}\\[0.3cm]
&&\displaystyle+2C_{YJ}\|\g\circ M\|_{BMO}\left\{E\left[\langle
|Z|\circ M\rangle_\sigma^T\, \biggm |\cF_\sigma
\right]\right\}^{1/2}. \ea\ee Using the elementary Cauchy
inequality, we have \be\ba{rcl} &&\displaystyle E\left[ \langle
Z\circ M+
M^\bot\rangle_\sigma^T|\cF_\sigma\right] \\[0.3cm]
&\le &\displaystyle
2\|\xi\|_{L^\infty}^2+2\|J\|_{\cR^\infty}^2+2C_{YJ}\|\a \circ
N_1\|_{BMO} \|
N_2\|_{BMO}\\[0.3cm]
&&\displaystyle +2C_{YJ}\|Y\|_{(\cR^\infty)^n}\|
\sqrt{\b}\circ M\|^2_{BMO}\\[0.3cm]
&&\displaystyle+2C_{YJ}^2\|\g\circ
M\|^2_{BMO}+{1\over2}E\left[\langle |Z|\circ M\rangle_\sigma^T\,
\biggm |\cF_\sigma \right]. \ea\ee The last inequality yields the
following \be\ba{rcl} &&\displaystyle {1\over 2}E\left[ \langle
Z\circ M+ M^\bot\rangle_\sigma^T\,
|\cF_\sigma\right]\\[0.3cm]
&\le &\displaystyle
2\|\xi\|_{L^\infty}^2+2\|J\|_{\cR^\infty}^2+2C_{YJ}\|\a \circ
N_1\|_{BMO} \|
N_2\|_{BMO}\\[0.3cm]
&&\displaystyle +2C_{YJ}\|Y\|_{(\cR^\infty)^n}\| \sqrt{\b}\circ
M\|^2_{BMO}+2C_{YJ}^2\|\g\circ M\|^2_{BMO}. \ea\ee Let $K$ denote
the right hand side of the last inequality.  We then have $Z\circ
M+ M^\bot\in BMO$ with $\|Z\circ M+ M^\bot\|_{BMO}\le \sqrt{2K}$.

Step 3. It remains to prove that $Z\circ M+ M^\bot\in
\overline{L^\infty}^{BMO}$. In view of the probabilistic version
of the Garnett and Jones theorem~\cite{GarnettJones} (due to
Varopoulos~\cite{Varo} and Emery~\cite{Emery3}, see
Kazamaki~\cite[Theorem 2.8, page 39]{Kaz}),  it is sufficient to
show that for any $\lambda>0$, \be
\sup_{\sigma}\left\|E\left[\exp{\left(\lambda |(Z\circ M)_T+
(M^\bot)_T-(Z\circ M)_\sigma- (M^\bot)_\sigma|\right)}|\,
|\cF_\sigma\right]\right\|_{L^\infty} <\infty.\ee Since
\be\ba{rcl} \displaystyle  Y_\sigma&=& \displaystyle \xi
+(J_T-J_\sigma)+\int_\sigma^Tf(s, Y_s)\, d\langle N_1, N_2
\rangle_s\\[0.3cm]
&& \displaystyle +\int_\sigma^T g(s,Y_s,Z_s)\, d\langle
M\rangle_s-\int_\sigma^TZ_s\, dM_s-\int_\sigma^T\, dM^\bot_s\ea
\ee and the random variable $Y_\sigma+J_\sigma-\xi-J_T\in
(L^\infty(\cF_T))^n$, it is sufficient to prove the following
\be\label{GJ} \sup_{\sigma}\biggm\|E\left[\exp{\left(\lambda
\biggm| -\int_\sigma^Tf(s, Y_s)\, d\langle N_1, N_2 \rangle_s
-\int_\sigma^T g(s,Y_s,Z_s)\, d\langle
M\rangle_s\biggm|\right)}\biggm|\cF_\sigma\right]\biggm\|_{L^\infty}<\infty.
\ee the left hand side of inequality~(\ref{GJ}) is equal to the
following \be\ba{rcl} & & \displaystyle
\sup_{\sigma}\biggm|E\left[\exp{\left(\lambda
\left(\int_\sigma^T|f(s, Y_s)|\, |d\langle N_1, N_2 \rangle_s|
+\int_\sigma^T |g(s,Y_s,Z_s)|\,
d\langle M\rangle_s\right)\right)}\biggm|\cF_\sigma\right]\biggm|_{L^\infty}\\[0.3cm]
&\le & \displaystyle \sup_{\sigma}\biggm|E\left[\exp{\left(\lambda
\int_\sigma^T \a_s\, |d\langle N_1, N_2 \rangle_s|+
\lambda\int_\sigma^T(\b_s |Y_s|+\g_s |Z_s|)\, d\langle
M\rangle_s\right)}\biggm|\cF_\sigma\right]\biggm|_{L^\infty}.
 \ea\ee
 While for any $\epsilon>0$
 \be\ba{rcl}&&\displaystyle  \lambda
\int_\sigma^T \a_s\, |d\langle N_1, N_2 \rangle_s|+
\lambda\int_\sigma^T(\b_s |Y_s|+\g_s |Z_s|)\, d\langle
M\rangle_s\\[0.3cm]
&\le & 2\epsilon \langle \a \circ N_1\rangle_\sigma^T
+2\epsilon^{-1}\lambda^2\langle N_2\rangle_\sigma^T
+\lambda\|Y\|_{(\cR^\infty)^n}\langle \sqrt{\b}\circ
M\rangle_\sigma^T\\[0.3cm]
&&\displaystyle +2\epsilon^{-1}\lambda^2\langle \g\circ
M\rangle_\sigma^T+2\epsilon \langle Z\circ M\rangle_\sigma^T,
\ea\ee in view of the facts that $b(N_2)=b(\sqrt{\b}\circ M)=b
(\g\circ M)=\infty$ (due to the assumption (ii) of the theorem),
 it is sufficient to prove the following for some $\epsilon>0$
\be \label{Est5-JN}\sup_{\sigma}\biggm|E\left[\exp{\left(
 4\epsilon \langle \a \circ N_1\rangle_\sigma^T+4\epsilon \langle Z\circ M\rangle_\sigma^T\right)}
 \biggm |\cF_\sigma\right]\biggm|_{L^\infty}<\infty.\ee
 Since $\a\circ M, Z\circ M\in BMO$, in view of the John-Nirenberg
 inequality (see Kazamaki~\cite[Theorem 2.2, page 29]{Kaz}), we
 have
 \be \sup_{\sigma}\biggm|E\left[\exp{\left(8\epsilon \langle \a \circ N_1\rangle_\sigma^T\right)}\biggm |\cF_\sigma\right]\biggm|_{L^\infty}
 \le {1\over {1-8\varepsilon \|\a\circ N_1\|^2_{BMO}}}<\infty\ee
 and
\be \sup_{\sigma}\biggm|E\left[\exp{\left(8\epsilon \langle Z\circ
M\rangle_\sigma^T\right)}\biggm
|\cF_\sigma\right]\biggm|_{L^\infty}
 \le {1\over {1-8\varepsilon \|Z\circ M\|^2_{BMO}}}<\infty\ee
 for sufficiently small $\varepsilon>0$. Therefore, the inequality~(\ref{Est5-JN}) hold when $\epsilon$ is sufficiently small. The proof is then complete.
\endpf

When the generator of a BSDE is not Lipschitz in the second unkown
variable, we should not expect that $Z\circ M\in
\overline{{L^\infty}}^{BMO}$  as in the last theorem. From
Kazamaki~\cite[Theorem 2.14, page 48]{Kaz}, we have
$\overline{{\cH^\infty}}^{BMO}\subset
\overline{{L^\infty}}^{BMO}$.  Therefore, we should not expect
that $Z\circ M\in \overline{{\cH^\infty}}^{BMO}$, neither. We have
the following negative result.

\bt Let $M$ be a one-dimensional standard Brownian motion, and
$\{\cF_t, 0\le t\le 1\}$ be the completed natural filtration.
Assume that $(Y, Z)$ solves the following BSDE:
\be\label{quadratic}\ba{rcl} dY_t&=&\displaystyle
Z_t\, dM_t+a Z_t^2\, dt, \quad t\in [0,1]; \\
Y_1&=&\xi \in L^\infty(\cF_1).\ea\ee  Then, $Y\in \cL^\infty$ and
$Z\circ M\in BMO$, but it is not always true that $Z\circ M\in
\overline{{L^\infty}}^{BMO}$.\et

{\bf Proof. } Without loss of generality, we assume $a={1\over
2}$. From Kobylansky~\cite{Kob1,Kob2} and Briand and
Hu~\cite{BrHu1}, we see that $Y\in \cL^\infty$ and $Z\circ M\in
BMO$.

Consider the  following process $X:$ \be X_t:=\int_0^t{1\over
\sqrt{1-s}}\, dM_s, \quad t\in[0,1).\ee Define the following
stopping time $\tau$: \be  \tau:=\inf\{\ t\in [0,1):
|X_t|^2>1\}.\ee It is easy to see that $\tau$ is a.s. well-defined
and $\tau<1$. Set \be \ba{rcl} \displaystyle \xi&:=&
-\log(X_\tau+2); \\
\displaystyle Y_t&:=& \displaystyle -\chi_{[0,\tau]}(t)\log
(X_t+2)+\xi \chi_{(\tau, 1]}(t), \quad
Z_t:=-{\chi_{[0,\tau]}(t)\over (X_t+2)\sqrt{1-t}}, \quad t\in
[0,1].\ea\ee Then, we can verify that $(Y,Z)$ is the unique
adapted solution of BSDE~(\ref{quadratic}). Further, in view of
the fact that $X_t+2\in [1,3]$, we have  \be
E\left[\exp{\left(\lambda \int_0^1Z_s^2\, ds\right)}\right]\ge
E\left[\exp{\left(\lambda \int_0^\tau{1\over 9(1-s)}\,
ds\right)}\right]=E\left[\exp{\left({\lambda \over 9}\langle X
\rangle_0^\tau \right)}\right]. \ee It is known (see
Kazamaki~\cite[Lemma 1.3, pages 11--12]{Kaz} for a similar result)
that \be E\left[\exp{\left({\lambda \over 9}\langle X
\rangle_0^\tau \right)}\right]=\infty \ee for $\lambda \ge {9\over
8}\pi^2.$ Consequently, we have \be E\left[\exp{\left(\lambda
\int_0^1Z_s^2\, ds\right)}\right]=\infty\ee for $\lambda \ge
{9\over 8}\pi^2$. In view of Kazamaki~\cite[Lemma 2.6, page
48]{Kaz} and BSDE~(\ref{quadratic}), we have $Z\circ M\not\in
\overline{{\cH^\infty}}^{BMO}$, and for $\lambda \ge {9\over
8}\pi^2$, \be E\left[\exp{\left(\lambda |Z\circ
W|\right)}\right]=\infty\ee
 due to both facts that $Y\in \cL^\infty$ and $\xi \in
 L^\infty(\cF_1)$.

Again, using the probabilistic version of the Garnett and Jones
theorem~\cite{GarnettJones} (see also Kazamaki~\cite[Theorem 2.8,
page 39]{Kaz}), we conclude the proof.
\endpf

\section{The linear case}

The study of linear BSDEs goes back to J. M. Bismut's Ph. D.
Thesis, which presented a  rather extensive study on stochastic
control,  optimal stopping, and stochastic differential games.
Also there, the concept of BSDEs was introduced and the theory of
linear BSDEs was initiated, though only for the case of uniformly
bounded coefficients and  $L^2$-integrable adapted solutions.

\subsection{BSDEs}

Assume that $A: \Omega\times [0,T]\to R^{n\times n}$ is $\{\cF_t,
0\le t\le T\}$-optional. Let $M$ be a continuous local martingale
such that $A\circ M\in BMO$. Consider the following linear SDE:
\be dX_t=A_tX_t\, dM_t+dV_t, \quad x_0=0. \ee

\bde \label{lsde}Consider the homogeneous linear SDE: \be
dX_t=A_tX_t\, dM_t, \quad X_0=I_{n\times n}. \ee Its unique strong
solution is denoted by $S(\cdot)$. It is said that $S(\cdot)$
satisfies the reversed H\"older inequality $(R_p)$ for some $p\in
[1,\infty)$ if for any stopping time $\sigma$ and any matrix norm
$|\cdot|$, we have
 \be E[|S(T)|^p|\cF_\sigma]\le C
|S(\sigma)|^p. \ee \ede

\br Note that $S=\cE(A\circ M)$ if $n=1$. In this case, it is
known that $S(\cdot)$ satisfies the reverse H\"older inequality
$(R_p)$ for all $p\in [1,\infty)$ if $A\circ M\in
\overline{L^\infty}^{BMO}$. See Kazamaki~\cite[Theorem 3.8, page
66]{Kaz} for details. Since $\langle B^T\rangle_T=T$ and thus
$B^T\in \cH^\infty\subset \overline{L^\infty}^{BMO}$, an immediate
consequence is the obvious fact that the stochastic exponential
$\cE (B^T)$ of a one-dimensional Brownian motion, stopped at a
deterministic time $T$, satisfies the reverse H\"older inequality
$(R_p)$ for all $p\in [1,\infty)$, which can be verified by some
straightforward explicit computations.\er

\br Assume that $A\circ M\in \overline{\cH^\infty}^{BMO}$. From
Theorem~\ref{nse}, we see that $S(\cdot)$ is uniformly integrable.
\er

Similar to the proof of Kazamaki~\cite[Corollary 3.2, page
60]{Kaz}, we can prove (by taking $U=|S(T)|^p$) the following
result.

\bt\label{continuation} Assume that $S(\cdot)$ is a uniformly
integrable matrix martingale, and let $p\in (1, \infty)$. If
$S(\cdot)$ satisfies the reverse H\"older inequality $(R_p)$, then
it satisfies $(R_{p'})$ for some $p'>p$. \et

We have the following \be
dS(t)^{-1}=-S(t)^{-1}[A_tdM_t-A_t^2d\langle M\rangle_t].\ee

\bt\label{lbsde} Let  $A\circ M\in BMO$ and $S(\cdot)$ be an
adapted continuous process that satisfies the reverse H\"{o}lder
inequality $(R_{p'})$ for $p'>1$. Let $q$ be the conjugate of
$p\in (1, p')$, i.e. ${1\over p}+{1\over q}=1$. Then, for $\xi\in
L^q(\cF_T)$ and $f\in \cL^q(0,T)$, the following BSDE
\be\label{lbsde eqn} \ba{rcl}
dY_t&=& -\left [A^\tau_t Z_t\, d\langle M \rangle_t+f_t\, dt\right]+Z_t\, dM_t+M_t^\bot, \quad \langle M, M^\bot\rangle=0,\\
Y_T&=&\xi \ea \ee has a unique adapted solution $(Y,Z\circ
M,M^\bot)\in (\cR^q)^n\times (\cH^q)^{2n}$.  Moreover,we have some
universal constant $K_q$ such that \be \label{aprioriestimate}
\|Y\|_{(\cR^q)^n}+\left\|\langle Y\rangle_T^{1/2}\right\|_{L^q}\le
K_q\left[~\|\xi\|_{(L^q)^n}+\|f\|_{(\cL^q)^n}\right]. \ee \et

\br In view of Theorem~\ref{continuation}, we can take $p=p'$ in
Theorem~\ref{lbsde} if furthermore $S(\cdot)$ is  assumed to be
uniformly integrable. \er

{\bf Proof of Theorem~\ref{lbsde}. } First for $s\ge t$, set \be
\widetilde Y_t:= E\left[S^\tau (t)^{-1}S^\tau (T)\xi
+\int_t^TS^\tau (t)^{-1}S^\tau (s) f_s\, ds
\biggm|\cF_t\right].\ee We have \be \widetilde Y_T=\xi \ee and \be
\widetilde Y_t=E\left[S^\tau (t)^{-1}S^\tau (T)\xi +\int_t^TS^\tau
(t)^{-1}S^\tau (T) f_s\, ds \biggm|\cF_t\right].\ee Since
$S(\cdot)$ satisfies the reverse H\"older inequality $(R_{p'})$,
letting $q'$ be the conjugate of $p'$, we see that \be\ba{rcl}
|\widetilde Y_t|&\le&\displaystyle E\left[|\xi|^{q'}\Bigm
|\cF_t\right]^{1/{q'}}
E\left [~\left|S^\tau (t)^{-1}S^\tau (T)\right |^{p'}\Bigm
|\cF_t\right]^{1/{p'}}\\[0.3cm]
&&\displaystyle +E\left[\int_t^T|f_s|^{q'}\,
ds\biggm|\cF_t\right]^{1/{q'}}E\left[\left|S^\tau (t)^{-1}S^\tau (T)\right|^{p'}\Bigm |\cF_t\right]^{1/{p'}}\\[0.3cm]
&\le&\displaystyle C\left(E\left[|\xi|^{q'}\Bigm
|\cF_t\right]^{1/{q'}}+E\left[\int_t^T|f_s|^{q'}\, ds\Bigm
|\cF_t\right]^{1/{q'}}\right).\ea\ee Therefore,  and using Doob's
inequality, we have \be\label{maximal estimate}\ba{rcl}
\left(E\left|\widetilde Y_t^*\right|^q\right)^{1/q}
&\le&\displaystyle
C\left(E\left[|\xi|^q\right]^{1/q}+E\left[\int_t^T|f_s|^q\,
ds\right]^{1/q}\right).\ea\ee Now it is clear that $\widetilde
Y\in (\cR^q)^n.$

We have \be \ba{rcl} d S^\tau (t)&=&S^\tau (t)A^\tau _t\, dM_t\\
dS^\tau (t)^{-1}&=&-\left[A^\tau _t\, dM_t-\left(A_t^{\tau
}\right)^2\, d\langle M\rangle_t\right]S^\tau (t)^{-1} \ea\ee and
\be S^\tau (t)\widetilde Y_t:= E\left[S^\tau (T)\xi +\int_0^T
S^\tau (s) f_s\, ds \biggm |\cF_t\right]-\int_0^tS^\tau (s)f_s\,
ds .\ee

From the martingale decomposition theorem, there is an $\{\cF_t,
0\le t\le T\}$-adapted process $z$ and a martingale $m^\bot$  such
that \be S^\tau (T)\xi +\int_0^T S^\tau (s) f_s\, ds=E\left[S^\tau
(T)\xi +\int_0^T S^\tau (s) f_s\, ds\right]+\int_0^Tz_s\,
dM_s+\int_0^Tdm^\bot_s, \quad \langle M, m^\bot\rangle=0. \ee
Then, we have \be S^\tau (t)\widetilde Y_t=E\left[S^\tau (T)\xi
+\int_0^T S^\tau (s) f_s\, ds\right]+\int_0^tz_s\,
dM_s+\int_0^tdm^\bot_s -\int_0^tS^\tau (s)f_s\, ds .\ee Denote  by
$X_t$ the right hand side of the last equality. Then, we have \be
\widetilde Y_t=S^\tau (t)^{-1}X_t, \quad dX_t=-S^\tau (t)f_t\,
dt+z_t\, dM_t+dm_t^\bot;\ee and from It\^o's formula, we further
have \be d\widetilde Y_t=-A_t^\tau Z_t\, d\langle M\rangle_t-f_t\,
dt +Z_t\, dM_t-dM_t^\bot\ee where \be Z_t:=A_t^\tau \widetilde
Y_t-S^\tau (t)^{-1}z_t, \quad M_t^\bot:=\int_0^tS^\tau
(s)^{-1}dm_s^\bot.\ee Noting that $\int_0^\cdot A_s^\tau Z_s\,
d\langle M\rangle_s$ is the quadratic variation of $\widetilde Y$
and the BMO martingale $A\circ M$, and then applying the a priori
estimate of Yor~\cite[Proposition 2, page 116]{Yor}, we have \be
E\left[\langle \widetilde Y\rangle^{q/2}\right]\le C_p
\left(1+\|A\circ M\|^q_{BMO}\right)E\left(|\widetilde
Y_T^*|^q\right). \ee The last inequality, together with
inequality~(\ref{maximal estimate}), shows that \be Z\circ M,
M^\bot\in (\cH^q)^n,\ee and the desired
estimate~(\ref{aprioriestimate}).  The proof for the existence is
complete.

The uniqueness follows immediately from the a priori
estimate~(\ref{aprioriestimate}).
\endpf

For the special case of $p=1$ (i.e, the conjugate number
$q=\infty$) and $n=1$, we have the following deeper result.

\bt\label{lbsdeinfty} Let  $A\circ M\in BMO$ and
$S(\cdot):=\cE(A\circ M)$ be its stochastic exponent.  Then, for
$\xi\in BMO(P)$ and $\int_0^T|f_s|\, ds\in BMO(P)$,
BSDE~(\ref{lbsde eqn}) has a unique adapted solution $(Y,Z\circ
M,M^\bot)\in (\cR^q)^n\times (\cH^q)^{2n}$ for any $q>1$.
Moreover, we have $Z\circ M+M^\bot\in BMO(P)$, and the following
estimate: \be \label{BMOaprioriestimate} \left\|Z\circ
M+M^\bot\right\|_{BMO}\le
C\left(\|\xi\|_{BMO}+\left\|\int_0^T|f_s|\,
ds\right\|_{BMO}\right) \ee for some universal constant $C$ which
depends on the BMO norm of $A\circ M$. \et

{\bf Proof of Theorem~\ref{lbsdeinfty}.}  The first assertion
follows immediately from Theorem~\ref{lbsde} and the fact that
$S(\cdot)$ satisfies the reverse H\"older inequality $(R_p)$ for
$p \in [1, p']$ for some $p'>1$. It remains to show the second
assertion. Without loss of generality, assume $f\equiv 0$.

First it is well known (see Kazamaki~\cite{Kaz}) that
$\xi^Q:\xi-\langle A\circ M, \xi\rangle\in BMO(Q)$ and $M^Q:
=M-\langle A\circ M, M \rangle\in BMO(Q)$ due to the fact that
$\xi, M\in BMO(P)$. In fact, for some $p>1$, $\cE(A\circ M)$
satisfies the reverse H\"older inequality $(R_p)$. See
Kazamaki~\cite{Kaz} for this assertion. Therefore, we have \be
\ba{rcl} &&\displaystyle  E_Q\left [\langle \xi^Q\rangle_t^T\Bigm |\cF_t
\right]\\[0.3cm]
&\le &\displaystyle E_Q\left [\langle \xi \rangle_t^T\Bigm |\cF_t \right] \\[0.3cm]
&=&\displaystyle  E\left[{\cE(A\circ M)_T\over \cE(A\circ M)_t}\langle \xi\rangle_t^T\biggm |\cF_t\right]\\[0.3cm]
&\le&\displaystyle  E\left[\biggm|{\cE(A\circ M)_T\over \cE(A\circ
M)_t}\biggm|^p\biggm |\cF_t\right]^{1/p}E\left[~\left|\langle
\xi\rangle_t^T\right|^q|\cF_t\right]^{1/q}\\[0.3cm]
&\le& \displaystyle C \|\xi\|_{BMO}, \ea \ee
 and the same is true for $M^Q$.

From BSDE~(\ref{lbsde eqn}), we have $\xi^Q=Z\circ M^Q+M^\bot\in
BMO(Q)$. Therefore, we have $Z\circ M+M^\bot\in BMO(P)$. Moreover,
we have the following estimate \be \left\|Z\circ
M+M^\bot\right\|_{BMO}\le C_1\|\xi^Q\|_{BMO(Q)}\le C_2
\|\xi\|_{BMO}. \ee The proof is then complete. \endpf

\subsection{SDEs}

For the multidimensional linear case, we have

 \bt\label{lsde} Let $A\circ M\in BMO$ such that $S(\cdot)$ satisfies the
reverse H\"{o}lder inequality $(R_{p'})$ for some $p'\in (1,
\infty)$, and $p\in (1,p')$. Then for $V\in (\cH^p)^n$, the
process \be X_t=S(t)\int_0^tS(s)^{-1}dV_s^Q, \quad 0\le t\le T\ee
solves SDE~(\ref{lsde}) and lies in $\cH^p$. Here, $V^Q:=V-\langle
A\circ M,V\rangle$. \et

\br In view of Theorem~\ref{continuation}, we can take $p=p'$ in
Theorem~\ref{lsde} if furthermore $S(\cdot)$ is  assumed to be
uniformly integrable. \er

The proof of Theorem~\ref{lsde} is based on a duality argument,
and will appeal to Theorem~\ref{lbsde}.

{\bf Proof of Theorem~\ref{lsde}.} Set \be \widetilde
X_t:=S(t)\int_0^tS(s)^{-1}dV_s^Q, \quad t\in [0,T]. \ee Using
It\^o's formula, we can show that \be d \widetilde
X_t=A_t\widetilde X_t dM_t+dV_t. \ee

For any $\xi\in (L^q(\cF_T))^n$ and $f\equiv 0$, BSDE~\ref{lbsde
eqn} has a unique adapted solution $(Y,Z,M^\bot)$. From It\^o's
formula, we have \be d (\widetilde X_t^\tau Y_t)=Y_t^\tau
dV_t+\widetilde X_t^\tau (Z_t\, dM_t+dM_t^\bot)+d\langle
Y,V\rangle_t, 0\le t\le T. \ee Therefore, applying the
inequality~(\ref{KWH}), we have \be E\langle \xi^\tau, \widetilde
X \rangle=E(\xi^\tau \widetilde X_T)=E\langle Y, V\rangle_T\le
\left\|\langle
Y\rangle_T^{1/2}\right\|_{L^q(\cF_T)}\|V\|_{(\cH^p)^n}.\ee In view
of the a priori estimate~(\ref{aprioriestimate}) of
Theorem~\ref{lbsde}, this shows $\widetilde X\in (\cH^p)^n.$\endpf

 As in Delbaen et al.~\cite{DMSSS1,DMSSS2},  we have the following theorem.

\bt\label{RHI} Let $p\in [1,\infty)$ and $A\circ M\in BMO$.
Suppose that the solution operator for SDE~(\ref{lsde}): \be
V\mapsto X=S(\cdot)\int_0^\cdot S(s)^{-1}dV_s^Q\in
\left(\cH^p\right)^n\ee is continuous from $(\cH^p)^n\to
(\cH^p)^n$. Here $V^Q:=V-\langle A\circ
M,V\rangle=(V_1,\cdots,V_n)^\tau-(\langle A\circ M,V_1 \rangle,
\cdots,\langle A\circ M, V_n \rangle)^\tau$ for
$V:=(V_1,\cdots,V_n)^\tau$. Then $S(\cdot)$ satisfies $(R_p(P))$.
Moreover, if $S(\cdot)$ is a uniformly integrable matrix
martingale, then the above solution operator remains to be
continuous from $\cH^{p'}\to \cH^{p'}$ for some $p'>p$.\et

{\bf Proof of Theorem~\ref{RHI}. } In view of
Theorem~\ref{continuation}, the second assertion is an immediate
consequence of the first one. Therefore, it is sufficient to prove
the first assertion.

For any stopping time $\sigma$, we are to show that \be
E\left[\left|S(T)S(\sigma)^{-1}\right|^p\, |\cF_\sigma\right]\le
C\ee for some constant $C$. For any $B\in \cF_\sigma$, take \be
V_i=\chi_{{[\sigma}, T]}\chi_B A\circ M, \quad i=1,2,\cdots,n.\ee
We have \be\ba {rcl} \|V\|_{(\cH^p)^n}^p&=&
E\left[\left(\langle A\circ M\rangle_T-\langle A\circ M\rangle_{\sigma}\right)^{p/2}\chi_B\right]\\[0.3cm]
&=& E\left[\chi_B E\left[\left(\langle A\circ M\rangle_T-\langle A\circ M\rangle_{\sigma}\right)^{p/2}\, \Bigm |\cF_{\widetilde \sigma}\right]\right]\\[0.3cm]
&\le & C\|A\circ M\|^p_{BMO}P(B).\ea\ee While \be\ba{rcl}
V^Q&=&\displaystyle \chi_{[{\sigma},T]}\chi_B A\circ M
-\chi_{[{\sigma},T]}\chi_B\langle A\circ M\rangle
(1,\cdots,1)^\tau,\\[0.3cm]
X_T&=&\displaystyle
S(T)\int_0^T\chi_{[{\sigma},T]}(s)S(s)^{-1}d(A\circ M)^Q\\
&=&\displaystyle -
S(T)\int_0^T\chi_{[{\sigma},T]}(s)\, d S(s)^{-1}\\
&=&\displaystyle - \chi_B S(T)\int_\sigma^T\, dS(s)^{-1}  \\[0.3cm]  &=&
S(T)\left[S(T)^{-1}-S({\sigma})^{-1}\right]\chi_B\\[0.3cm]
&=&\displaystyle  \chi_B\left[S(T)S({\sigma})^{-1}-I\right]. \ea
\ee From the assumption of the underlying theorem, we have \be
\left\| \chi_B\left[S(T)S({
\sigma})^{-1}-I\right]\right\|_{(\cH^p)^{n\times
n}}=\|X_T\|_{(\cH^p)^{n\times n}}\le C \|V\|_{(\cH^p)^n}\le
C\|A\circ M\|_{BMO_p}P(B) \ee for some constant $C>0.$ Therefore,
in view of the BDG inequality,  the quantity $$\left\|S(T)S({
\sigma})^{-1}\chi_B\right\|_{(L^p)^{n\times n}}$$ is bounded by
$P(B)$. This implies that $S(\cdot)$ satisfies the reverse
H\"older inequality $(R_p)$.
\endpf

\br From Theorems~\ref{nsesliceable} and~\ref{RHI}, we see that if
$A\circ M\in \overline{\cH^\infty}^{BMO}$, then $S(\cdot)$
satisfies the reverse H\"older inequality $(R_p)$ for all $p\in
[1, \infty)$.\er


For the special case of $p=1$ and $n=1$,  we have the following
better result.

\bt\label{lsde1} Let $n=1$. Assume that $A\circ M\in BMO$. Then
for $V\in \cH^1$, the local martingale  \be
X_t=S(t)\int_0^tS(s)^{-1}dV_s^Q, \quad t\in [0,T]\ee solves
SDE~(\ref{lsde}) and lies in $\cH^1$. Here, $V^Q:=V-\langle A\circ
M,V\rangle$. \et

{\bf Proof of Theorem~\ref{lsde1}. }  Take any $\xi \in BMO$ and
$f\equiv 0$. Let $(Y, Z\circ M, M^\bot)$ be the unique solution of
BSDE~(\ref{lbsde eqn}) for the data $(\xi,f)$. As shown in the
proof of Theorem~\ref{lsde}, we have \be E\langle \xi^\tau,
\widetilde X \rangle=E(\xi^\tau \widetilde X_T)=E\langle Y,
V\rangle_T=E\langle Z\circ M+M^\bot , V\rangle_T.\ee Applying
Fefferman's inequality, we have \be E\langle \xi^\tau, \widetilde
X \rangle\le \sqrt{2}\left\|Z\circ
M+M^\bot\right\|_{BMO}\|V\|_{(\cH^1)^n}.\ee In view of the a
priori estimate~(\ref{BMOaprioriestimate}) of
Theorem~\ref{lbsdeinfty}, we have \be E\langle \xi^\tau,
\widetilde X \rangle\le C \|\xi\|_{BMO}\|V\|_{(\cH^1)^n}, \quad
\forall \xi\in BMO\ee for some positive constant $C$. In view of
Lemma~\ref{1-infty
 duality}, The last inequality implies that
 $\widetilde X\in (\cH^1)^n$. The proof is then complete.
\endpf

\section{One-dimensional linear case: the characterization of Kazamaki's critical quadratic exponent being infinite.}

In Section 2, we have applied Fefferman's inequality to prove new
results for SEs and BSDEs.  In what follows, we present an
operator approach to Kazamaki's critical quadratic exponent on BMO
martingales. We establish some relations between Kazamaki's
critical quadratic exponent $b(M)$ of a BMO martinagle $M$ and the
solution operator for the associated $M$-driven SDE. Throughout
this section, all processes will be considered in $[0, \infty).$

Let $M \in BMO$ be real and $p\in [1, \infty).$ Consider the
operator $\phi: \phi(X)=X\circ M$ for $X\in \cH^p$. Define the
complex version $\tilde \phi: \cH^p(\mathbb{C})\to
\cH^p(\mathbb{C})$ as follows: \be\tilde \phi (U+iV):=U\circ
M+iV\circ M.  \ee Since \be\ba{rcl}  \|\tilde \phi (U+iV)\|_p&\le&
\|\phi (U)\|_{\cH^p}+ \|\phi
(V)\|_{\cH^p}\\
&\le & \|\phi\|\|U\|_{\cH^p}+\|\phi\|\|V\|_{\cH^p}\le
2\|\phi\|\|U+iV\|_{\cH^p}, \ea \ee we have \be \|\tilde \phi\|\le
2\|\phi\|, \quad \|\tilde \phi^n\|\le 2\|\phi^n\|\quad \hbox {\rm
(since $M$ is real!)}.\ee Their spectral radii are equal, denoted
by $r_p$: \be \lim_{n\to \infty}\|\tilde \phi^n\|^{1/n}=\lim_{n\to
\infty}\|\phi^n\|^{1/n}=r_p.\ee

For $\l\in \mathbb{C}$, define \be M^\l: \l M-\l^2\langle
M\rangle\ee and \be \cE(\l M)_t:=\exp\left(\l M_t-{1\over 2}\l^2
\langle M\rangle_t\right)\ee which is a complex local martingale.

Using the same procedure as in the real case, we have

\bp\label{crh} Suppose that $\l\in \mathbb{C}$ and that $(Id-\l
\tilde \phi)$ has an inverse on $\cH^p(\mathbb{C})$ for some $p\in
[1, \infty)$. Then $\cE(\l M)$ satisfies the following stronger
property than ($R_p$): there is a positive constant $K$ such that
\be E\left[\left|{\cE(\l M)_\sigma\over \cE(\l M)_\tau}
\right|^p\, \biggm|\cF_\tau\right]\le K\ee for any stopping times
$\tau$ and $\sigma$ such that $0\le \tau\le \sigma\le \infty$. \ep

\br Even for $p=1$, Proposition \ref{crh} yields information in
the complex case.\er

{\bf Proof of Proposition \ref{crh}. } Use stopping to make all
integrals  bounded. For $A\in \cF_\tau$, define the process $g$ as
follows: \be g(t)=\chi_A\chi_{(\tau,\infty)}(t),\quad t\in
[0,\infty).\ee Then, we have \be \|g\circ M\|_{\cH^p}\le C
[P(A)]^{1\over p}.\ee

Indeed, we have \be \ba{rcl}E\left[\langle g\circ
M\rangle_\infty^{p/2}\right]&=&\displaystyle  E\left[\chi_A
\left(\int_\tau^\infty
d\langle M\rangle_t\right)^{p/2}\right]\\[0.5cm]
&=&\displaystyle    E\left[\chi_A (\langle
M\rangle_\tau^\infty)^{p\over 2}\right]
\\
&\le &\displaystyle   C P(A). \qquad (\hbox{\rm noting that $M\in$
BMO})\ea \ee

On the other hand, we have \be \ba{rcl}\displaystyle
\phi^\l(g\circ M)_\sigma &:=&\displaystyle \cE(\l M)_\sigma
\int_\tau^\sigma \cE(\l
M)_s^{-1}g(s) \, dM_s^\l\\[0.3cm]
&=&\displaystyle  \displaystyle \cE(\l
M)_\sigma \int_\tau^\sigma \cE(\l M)_s^{-1}\chi_A \, dM_s^\l\\[0.3cm]
&=&\displaystyle  \cE(\l M)_\sigma \left(\cE(-\l
M^\l)_\tau-\cE(-\l
M^\l)_\sigma\right) \chi_A \\[0.3cm]
&=&\displaystyle  \chi_A \left({\cE(\l M)_\sigma\over \cE(\l
M)_\tau}-1\right). \ea \ee Since the map $\phi^\l: g\circ M\to
\phi^\l(g\circ M)$ is the operator $(Id-\l \phi)$, we get by
hypothesis that there is a constant $K$ (changing from line to
line) such that \be E\left[\chi_A\left|{\cE(\l M)_\sigma \over
\cE(\l M)_\tau}-1\right|^p\right]\le K \|g\circ M\|_{\cH^p}^p\le K
P(A).\ee Therefore, \be E\left[\left|{\cE(\l M)_\sigma\over \cE(\l
M)_\tau}-1\right|^p\biggm | \cF_\tau\right]\le K, \ee which
implies the following \be E\left[\left|{\cE(\l M)_\sigma\over
\cE(\l M)_\tau}\right|^p\biggm | \cF_\tau\right]\le K. \ee
\endpf

\bp\label{crhr} Suppose that $\cE(\l M)$ satisfies ($R_p$) for
some $\l\in \mathbb{C}$ and some $p\in (1, \infty)$, that is,
there is a positive constant $K$ such that \be
E\left[\left|{\cE(\l M)_\sigma\over \cE(\l M)_\tau} \right|^p\,
\biggm |\cF_\tau\right]\le K\ee for any stopping times $\tau$ and
$\sigma$ such that $0\le \tau\le \sigma\le \infty$. Then $(Id-\l
\tilde \phi)$ has an inverse on $\cH^p(\mathbb{C})$.  \ep

Kazamaki~\cite[Lemma 2.6, page 48]{Kaz} states that
 \be {1\over
\sqrt{2} d_2(M, \cH^\infty)}\le b(M).\ee
Schachermayer~\cite{Schacher} has shown that the reverse is not
true in the following sense: \be b(M)=+\infty \nRightarrow M\in
\overline {\cH^\infty}^{BMO}.\ee There seems to be no hope to
establish a relation between $\hbox{\rm dist}(M, \cH^\infty)$ and
$b(M)$.

\bp \label{ubd} If $\l\in \mathbb{C}$ satisfies \be |\l|<
{b(M)\over \sqrt{2p(2p-1)}},\ee then $(Id-\l \tilde \phi)$ has an
inverse on $\cH^p(\mathbb{C})$. \ep

{\bf Proof of Proposition \ref{ubd}.} In view of
Proposition~\ref{crhr}, it is sufficient to show that $\cE(\l M)$
satisfies ($R_p$), i.e., there is a positive constant $C$ such
that \be \label{237} E\left[\left|\cE(\l N)\right|^p\biggm
|\cF_\tau\right]\le C\ee with $N:=M-M^\tau$.

Denote $\l:=u+iv$ with $u$ and $v$ being real numbers. We have
\be\ba{rcl}&&\displaystyle
\left|\exp{\left(p(u+iv)N_\infty-{1\over2}(u+iv)^2p\langle N
\rangle_\infty\right)}\right|\\
&=& \displaystyle
\exp{\left(puN_\infty-{1\over2}p(u^2-v^2)\langle N
\rangle_\infty\right)}\\
&=&\displaystyle \exp{\left(puN_\infty-p^2u^2\langle
N\rangle_\infty\right)} \exp{\left(p^2u^2\langle
N\rangle_\infty-{1\over2}pu^2\langle N
\rangle_\infty+{1\over2}pv^2\langle N \rangle_\infty\right)}.
\ea\ee Taking the conditional expectation and using the
Cauchy-Schwarz inequality, we have \be\ba{rcl}&& \displaystyle
E\left[\left|\cE(\l N)_\infty\right|^p\biggm |\cF_\tau\right]\\
&\le & \displaystyle E\left[\exp{\left(2pu N_\infty-2p^2u^2\langle
N
\rangle_\infty\right)}|\cF_\tau\right]^{1\over2}E\left[\exp{\left(\langle
N
\rangle_\infty\left(2p^2u^2-pu^2+pv^2\right)\right)}\right]^{1\over2}.\ea
\ee Since $M$ is in BMO and so is $2puM$, $\cE(2puM)$ is uniformly
integrable. Hence, we have \be E\left[\exp{\left(2pu
N_\infty-2p^2u^2\langle N
\rangle_\infty\right)}|\cF_\tau\right]=1.\ee Concluding the above,
we have \be E\left[\left|\cE(\l N)_\infty\right|^p\biggm
|\cF_\tau\right]\le E\left[\exp{\left(\langle N
\rangle_\infty\left(2p^2u^2-pu^2+pv^2\right)\right)}\right]^{1\over2}.
\ee In view of the fact that \be 2p^2u^2-pu^2+pv^2\le
p(2p-1)(u^2+v^2)=p(2p-1)|\l|^2< {1\over2}b^2(M),\ee we obtain the
desired inequality~(\ref{237}).
\endpf

The spectral radius $r_p$ of $\tilde \phi: \cH^p(\mathbb{C}) \to
\cH^p(\mathbb{C})$ is estimated by $b(M)$.

\bc We have \be r_p\le {\sqrt{2p(2p-1)}\over b(M)}.\ee\ec

{\bf Proof. } Since \be (1-\l \tilde \phi)^{-1} \hbox{ \rm exists
} \Longleftrightarrow (\l^{-1}-\tilde \phi)^{-1} \hbox{ \rm exists
},\ee we have from Proposition \ref{ubd} that \be |\l^{-1}|>r_p\ee
for all $\l\in \mathbb{C}$ such that \be |\l|< {b(M)\over
\sqrt{2p(2p-1)}}.\ee This implies immediately the desired
inequality.\endpf

\bp \label{lbd} We have \be r_p\ge {\sqrt{p}\over b(M)}. \ee\ep

{\bf Proof of Proposition \ref{lbd}.  }  Take $\l\in R$ such that
$\l< r_p^{-1}.$ Then $(Id-i\l \tilde \phi)$ has inverse. From
Proposition~\ref{crh}, we see that there is a positive constant
$C$ such that \be E\left[\left|{\cE(\l M)_\infty\over \cE(\l
M)_\tau} \right|^p\biggm |\cF_\tau \right]\le C.\ee
\endpf
Set $N:=M-M^\tau$. We have \be E\left[\left|\exp{\left(ip\l
N_\infty+{1\over2}p\l^2\langle N\rangle \right)}\right|\biggm
|\cF_T\right]\le C, \ee which is equivalent to the following \be
E\left[\exp{\left({1\over2}p\l^2\langle N\rangle \right)}\biggm
|\cF_T\right]\le C. \ee Therefore, by the definition of $b(M)$, we
have \be p\l^2\le b^2(M)\ee for all $\l\in R$ such that $\l<
r_p^{-1}.$ The desired inequality then follows immediately.
\endpf

We combine  the above two propositions into the following theorem

\bt We have \be {\sqrt{p}\over b(M)}\le r_p\le {\sqrt{2p
(2p-1)}\over b(M)}.\ee\et

We have the following equivalent conditions.

 \bt The following
statements are equivalent.

{\rm (i)} $\forall \l\in \mathbb{C}, \forall p\in [1, \infty)$,
the map $Id-\l \phi: \cH^p\to \cH^p$ is an isomorphism .

{\rm (ii)} For some $p\in [1, \infty)$, the map $Id-\l \phi:
\cH^p\to \cH^p$ is an isomorphism for $\forall \l\in \mathbb{C}$.

{\rm (iii)} $\forall \l\in \mathbb{C}, \forall p\in [1, \infty)$,
there is $K>0$ such that \be E\left[\left|{\cE(\l M)_\sigma\over
\cE(\l M)_\tau} \right|^p\biggm |\cF_\tau\right]\le K\ee
 for any  stopping times $\tau$ and
$\sigma$ such that $0\le \tau\le \sigma\le \infty$.

{\rm (iv)}  For some $p\in [1, \infty)$ and $\forall \l\in
\mathbb{C}$, there is $K>0$ such that \be E\left[\left|{\cE(\l
M)_\sigma\over \cE(\l M)_T} \right|^p\biggm |\cF_\tau\right]\le
K\ee for any stopping times $\tau$ and $\sigma$ such that $0\le
\tau\le \sigma\le \infty$.

{\rm (v)} $\tilde \phi: \cH^p(\mathbb{C})\to \cH^p(\mathbb{C})$ is
quasinilpotent (i.e. $r_p=0$) for $\forall p\in [1, \infty)$.

{\rm (vi)} $\tilde \phi: \cH^p(\mathbb{C})\to \cH^p(\mathbb{C})$
is quasinilpotent (i.e. $r_p=0$) for some $p\in [1, \infty)$.

{\rm (vii)} $b(M)=+\infty$.

{\rm (viii)} $\displaystyle \lim_{n\to \infty} \|\tilde
\phi^n\|^{1/n}=\lim_{n\to \infty}\|\phi^n\|^{1/n}=0.$\et

{\bf Proof. } We show that (vii) $\Longrightarrow$ (iii). For
$\forall \l=\l_1+i\l_2\in \mathbb{C}$, we have
\be\ba{rcl}\displaystyle  \left|{\cE(\l M)_\infty\over \cE(\l
M)_\tau}\right|&=&\displaystyle
\left|\exp{\left[\l(M_\infty-M_\tau)-{1\over
2}(\langle M\rangle_\infty-\langle M\rangle _\tau)\right]}\right|
\\[0.3cm]
&=&\displaystyle \exp{\left[\l_1(M_\infty-M_\tau)-{1\over2}(\l_1^2-\l_2^2)(\langle M\rangle_\infty-\langle M\rangle _\tau)\right]}\\[0.3cm]
&\le&\displaystyle
\exp{[|\l_1||M_\infty-M_\tau|]}\exp{\left[{1\over2}\l_2^2(\langle
M\rangle_\infty-\langle M\rangle _\tau)\right]}.\ea\ee Hence,
\be\ba{rcl}\displaystyle  \left|{\cE(\l M)_\infty\over \cE(\l
M)_\tau}\right|^p &\le&\displaystyle
\exp{[|\l_1|p|M_\infty-M_\tau|]}\exp{\left[{1\over2}\l_2^2p(\langle
M\rangle_\infty-\langle M\rangle _\tau)\right]}.\ea\ee

On the other hand, we have for $\forall \l=\l_1+i\l_2\in
\mathbb{C}$, \be\ba{rcl}\displaystyle
\exp{\left[{1\over2}(\l_2^2-2\l_1^2)\langle M
\rangle^\infty_\tau\right]}&=&\displaystyle \left|\cE(\l
M^\infty_\tau)\right|\cE{(\l_1M^\infty_\tau)},\ea\ee from which we
can derive that (iii) $\Longrightarrow$ (vii).

\end{document}